\documentclass[11pt]{article}
\usepackage{graphicx, setspace, amsmath, amssymb, amsthm, amsfonts, url, amsthm, accents, gensymb, enumitem, mathtools, listings, color, fullpage, listings, courier, longtable}
\usepackage[linesnumbered,ruled,vlined]{algorithm2e}
\usepackage[english]{babel}
\usepackage[utf8x]{inputenc}
\newtheorem{theorem}{Theorem}
\newtheorem{definition}{Definition}

\renewenvironment{thebibliography}[1]{%
	\begin{oldthebibliography}{#1}%
		\setlength{\parskip}{0ex}%
		\setlength{\itemsep}{0ex}%
	}%
	{%
	\end{oldthebibliography}%
}
\definecolor{mygreen}{rgb}{0,0.6,0}
\definecolor{mygray}{rgb}{0.85,0.85,0.85}
\definecolor{mymauve}{rgb}{0.58,0,0.82}
\SetCommentSty{mycommfont}

\title{Graph Colouring: A Visual Tour}
\author{Rhyd Lewis\\
	School of Mathematics, Cardiff University, CF24 4AX, Wales. \\
	\url{LewisR9@cardiff.ac.uk}
}

\begin{document}
	
	\maketitle
	
	\begin{abstract}
		Graph colouring is a combinatorial optimisation problem with applications in several important domains, including sports scheduling, cartography, street map navigation, and timetabling. It is also of significant theoretical interest and a standard subject in university-level courses on graph theory, algorithms, and combinatorics. In this paper, we consider the topics of node, edge, and face colouring along with their associated algorithms. Theoretical results are reviewed and brought to life through a collection of detailed, visually engaging figures designed to enhance understanding and appeal.
	\end{abstract}
	
	
	\section{Introduction}\label{sec:intro}
	
	Graph colouring refers to the computational task of assigning colours to the elements of a graph so that pairs of adjacent elements never share the same colour. The aim is to minimise the number of colours that are used. 
	
	In broad terms, a graph is a mathematical structure comprising a collection of \emph{nodes} in which some pairs of nodes are connected by \emph{edges}. More formally, a graph $G$ is defined by an ordered pair $G=(V,E)$, where $V$ is a set of $n$ nodes, $E$ is a set of $m$ edges, and each edge represents a connection between a pair of nodes. Nodes are also known as \emph{vertices}. Two nodes $u,v\in V$ are described as \emph{adjacent} if and only if they are joined by an edge; that is, $\{u,v\}\in E$. The \emph{neighbours} of a node are all the nodes adjacent to it, and the degree of node $v$, denoted by $\deg(v)$, is the number of neighbours. Similarly, a pair of edges is adjacent if and only if they share a common endpoint, and the neighbours of an edge $\{u,v\}$ are all edges containing the endpoints $u$ or $v$.
	
	\begin{figure}[htbp]
		\begin{center}
			\includegraphics[trim= 10 10 10 10, clip, scale=0.44]{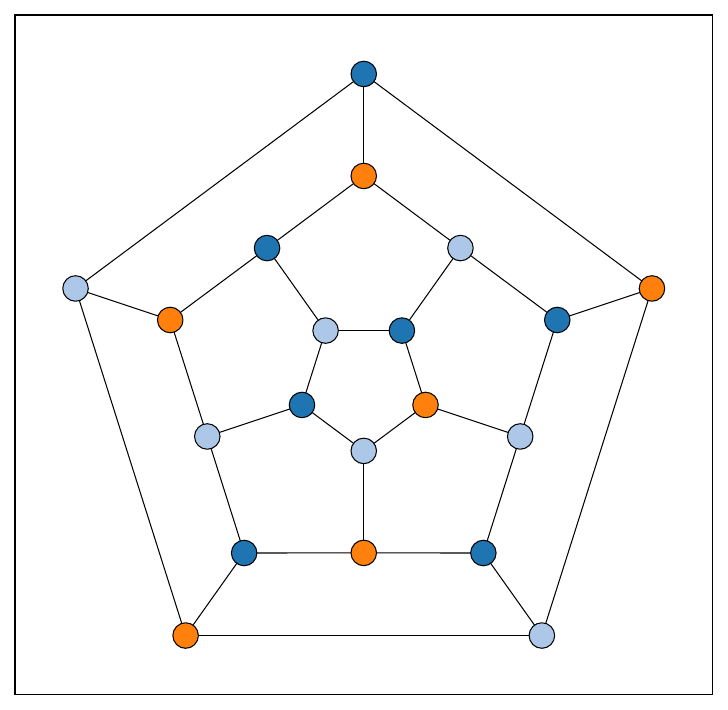}
			\includegraphics[trim= 10 10 10 10, clip, scale=0.44]{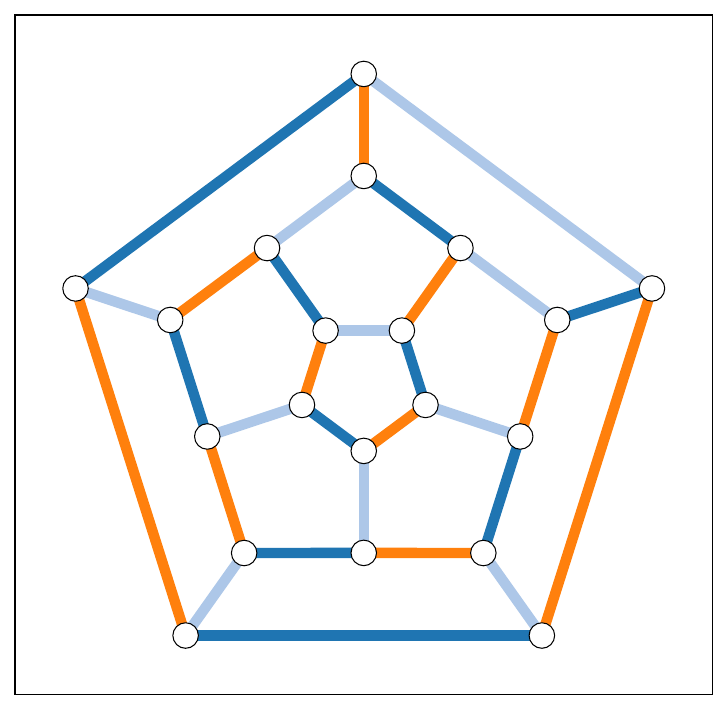}
			\includegraphics[trim= 10 10 10 10, clip, scale=0.44]{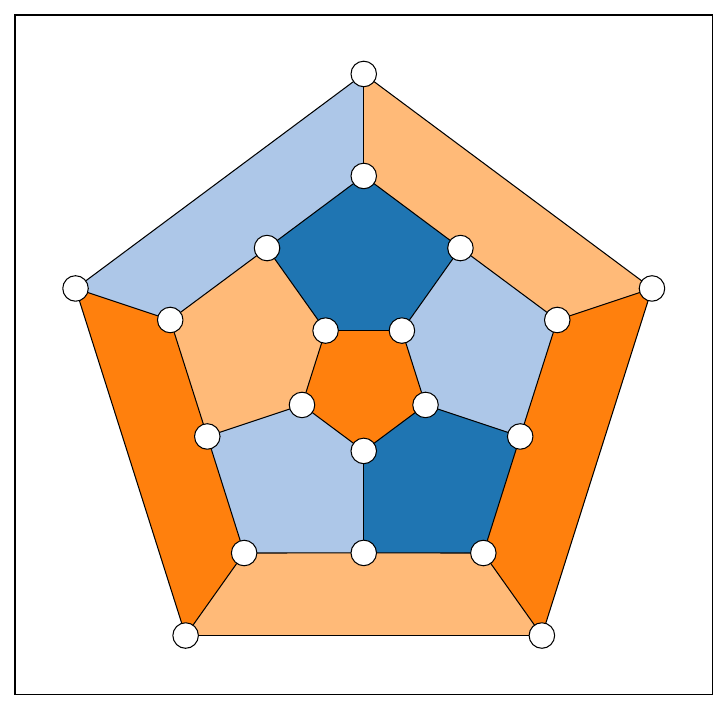}
		\end{center}
		\caption{A node colouring, edge colouring, and face colouring (respectively) of the dodecahedral graph. This graph has $n=20$ nodes and $m=30$ edges. As shown, its chromatic number and chromatic index are both three. Because it is planar, the faces of this graph can also be coloured. In this case, the face chromatic number is four.}
		\label{fig:dodec}
	\end{figure}
	
	Examples of these concepts are provided in Figure~\ref{fig:dodec}, which also demonstrate the following definitions.
	
	\begin{definition}
		\label{def:nodecol}
		A \emph{node colouring} of a graph $G$ is an assignment of colours to nodes so that all pairs of adjacent nodes have different colours. The smallest number of colours needed to colour the nodes of $G$ is known as its \emph{chromatic number}, denoted by $\chi(G)$.
	\end{definition}
	
	\begin{definition}
		\label{def:edgecol}
		An \emph{edge colouring} of $G$ is an assignment of colours to edges so that all pairs of adjacent edges have different colours. The smallest number of colours needed to colour the edges of $G$ is known as its \emph{chromatic index}, denoted by $\chi'(G)$. 
	\end{definition}
	
	\begin{definition}
		\label{def:facecol}
		A \emph{face colouring} of a graph is an assignment of colours to the faces of one of its planar embeddings (if such an embedding exists) so that faces with common boundaries have different colours. The smallest number of colours needed to colour the faces of a planar embedding of $G$ is known as its \emph{face chromatic number}, denoted by $\chi_f(G)$.
	\end{definition}
	
	Face colourings require nodes to be arranged on the plane so that none of the graph's edges intersect. Consequently, they are only possible for planar graphs. A colour should also be assigned to the region surrounding the graph, usually known as the unbounded face. For presentational purposes, the colour of the unbounded face is not shown in the rightmost graph of Figure~\ref{fig:dodec}; however, it can be assigned to light blue here, maintaining the use of four colours. Euler's formula also tells us that the number of faces $f$ in a planar embedding is strictly related to the number of nodes $n$ and edges $m$. Specifically, $n-m+f=2$. 
	
	Due to the famous four colour theorem, it is known that $\chi_f(G)$ never exceeds four~\cite{Appel1977}. The problem of identifying face colourings that use four colours is also solvable in quadratic time~\cite{Robertson1997}, though the problem of determining if $\chi_f(G)=3$ is $\mathcal{NP}$-complete in general~\cite{Dailey1980}. Similarly, the problems of identifying optimal node and edge colourings (that is, colourings that use the minimum number of colours) are known to be $\mathcal{NP}$-hard~\cite{Garey1979}.
	
	Node colouring can be considered the most fundamental of the graph colouring problems, because edge and face colouring problems can always be converted into instances of the former. Specifically:
	\begin{itemize}
		\item An edge colouring of a graph $G$ can be achieved by colouring the nodes of its \emph{line graph} $L(G)$, where each node of $L(G)$ corresponds to an edge of $G$, and two nodes in $L(G)$ are adjacent if and only if the corresponding edges in $G$ are adjacent. As such, the colour of each node in $L(G)$ gives the colour of the corresponding edge in $G$, implying that $\chi'(G)=\chi(L(G))$. 
		\item A face colouring of a planar graph $G$ can be found by colouring the nodes of its \emph{dual graph} $G^*$. To form $G^*$, a planar embedding of $G$ is taken. Each node in $G^*$ then corresponds to a face in the embedding, and pairs of nodes in $G^*$ are made adjacent if and only if the corresponding faces share a common boundary. The colours of the nodes in $G^*$ then give the colours of the faces in the embedding. Hence, $\chi_f(G)=\chi(G^*)\leq 4$.
	\end{itemize}
	
	Graph colouring is a topic of theoretical interest that is often taught in university-level courses on graph theory, algorithms, and combinatorics~\cite{Cranston2024,Lewis2021}; indeed, before being proved in the 1970s, the four colour theorem was one of the most famous unsolved problems in all of mathematics~\cite{Wilson2003}. Our aim in this paper is to review and visualise fundamental results from this field. We also want to present aesthetically pleasing figures arising from these results using, among other things, fractals, triangulations, and tiling patterns. To do this, we make extensive use of the open-source Python library \emph{GCol}, which is maintained by the author of this paper~\cite{Lewis2026GCol, Lewis2025}. This library provides several high-performance algorithms for graph colouring, together with tools for solution visualisation. Most graphs considered in this paper have been generated by the author; the remaining cases are taken from the House of Graphs database \cite{Brinkmann2013} and are indicated in the text by the prefix ``HoG'' followed by their identification number. Colours used in all figures make use of the Tableau palette, which is intended to provide a set of colours that are visually distinct, perceptually balanced, and colour-blind-friendly.
	
	In the following three sections, we consider the node, edge, and face colouring problems in turn. In Section~\ref{sec:nodecol}, we start by discussing the node colouring algorithms used in this work and examine their relative strengths and drawbacks. These algorithms are also used to produce the edge and face colourings given in Sections~\ref{sec:edgecol} and \ref{sec:facecol} by node colouring the graphs' line and dual graphs, respectively. The code and data used to produce all results and images can be found at \url{https://doi.org/10.5281/zenodo.18713617}. 
		
	\section{Node Colouring}
	\label{sec:nodecol}
	
	As stated previously, the $\mathcal{NP}$-hard node colouring problem requires solutions to use a minimum number of colours while ensuring that adjacent nodes have different colours. As such, the problem can also be viewed as the task of partitioning a graph's nodes into a minimum number of independent sets. Practical applications of node colouring have been noted in several areas, including timetabling~\cite{Cambazard2012, Lewis2008}, the design of seating plans~\cite{Lewis2016}, frequency assignment problems~\cite{Aardel2002, Dupont2009} and solving sudoku problems. 
	
	\subsection{Algorithms for Node Colouring}
	
	The 2021 book of Lewis~\cite{Lewis2021} surveys many different algorithms for node colouring, including exponential-time exact algorithms based on backtracking and integer programming, polynomial-time constructive heuristics, and metaheuristic-based approaches. Perhaps the most basic approach is the well-known greedy colouring algorithm. This assumes that colours are identified by the labels $\{0,1,2,\ldots\}$; the method then considers each node in turn and simply assigns it to the lowest colour label not being used by any of its neighbours. If the node ordering is determined before execution, the greedy method has a complexity of $O(n+m)$. It is also known that there is always at least one ordering that leads to an optimal node colouring, though identifying this ordering is also an $\mathcal{NP}$-hard problem~\cite{Lewis2021}.
	
	\begin{algorithm}[bt!]
		\label{alg:bktr}
		\DontPrintSemicolon
		\SetKwInOut{Input}{input}
		\SetKwInOut{Output}{output}
		\ResetInOut{output}
		\caption{Backtracking Algorithm for Node Colouring}
		\Input{A graph $G=(V,E)$.}
		\Output{An optimal node colouring $\mathcal{S}$.}
		Let $C=\{v_1,v_2,\ldots\}$ be a large clique in $G$ and let $k=n$.\;
		\ForEach{$v_i\in C$}{
			Assign node $v_i$ to colour $i$.\;
		}
		Choose an uncoloured node $u$ in $G$.\;
		\textsc{Colour}($u$)\;
		\Return $\mathcal{S}$
		
		\bigskip
		\underline{\textbf{subroutine} \textsc{Colour}($u$)}\;
		Let $l$ be the number of colours currently being used in the colouring of $G$.\;
		\If{\textup{all nodes in $G$ are coloured}}{
			Let $\mathcal{S}$ be a copy of the current node colouring.\;
			\lIf{$l=|C|$}{\Return True}
			\lElse{set $k=l-1$ and \Return False}
		}
		Set $i=0$\;
		\While{\textup{True}}{
			\If{$i \geq k$ \textup{\textbf{or}} $l > k$}{
				\textbf{break}\;
			}
			\ElseIf{\textup{colour $i$ is feasible for node $u$}}{
				Assign node $u$ to colour $i$.\;
				Choose an uncoloured node $v$ in $G$.\;
				\If{\textup{\textsc{Colour}($v$) = True}}{
					\Return True
				}
				Remove the current colour from node $u$.\;
			}
			Set $i=i+1$\;
		}
		\Return False
	\end{algorithm}
	
	The greedy colouring algorithm can also be extended to form an exponential-time recursive backtracking algorithm that guarantees an optimal solution. A specification of this is given in Algorithm~\ref{alg:bktr}. As shown, in each recursive step a currently uncoloured node $u$ is selected. Each available colour $i\in\{0,1,\ldots,k-1\}$ is then considered in turn and, for each feasible value of $i$, the algorithm tentatively assigns $u$ to colour $i$ before moving on to the next uncoloured node $v$. If no suitable colour exists for $u$, the algorithm backtracks by undoing the most recent assignment and trying the next colour. The process then continues, systematically trying all feasible colours with all nodes. During execution, the number of available colours $k$ is also decremented whenever a full $k$-colouring of $G$ is achieved. This ensures that, at the end of execution, $\chi(G)$ colours are being used in the returned node colouring $\mathcal{S}$. Line~1 of Algorithm~\ref{alg:bktr} also contains the optional step of identifying and colouring a clique $C$ in $G$.\footnote{Given the graph $G=(V,E)$, a clique $C$ is a subset of nodes such that, for all $u,v \in C$, $\{u,v\}\in E$. Calculating the \emph{largest} in $G$ is an $\mathcal{NP}$-hard problem, though smaller cliques can be calculated in polynomial time.} All nodes in a clique must be assigned to different colours. Consequently, the cardinality of $C$ provides a lower bound on $\chi(G)$, allowing the backtracking process to halt if this bound is achieved.
	
	In our descriptions of the greedy and backtracking algorithms above, the order in which nodes are selected to be coloured is left unspecified; however, it is known that some orderings can cause the greedy algorithm to use many additional colours and the backtracking algorithm to exhibit very poor run times~\cite{Johnson1974,Kubale1985}. An effective strategy for mitigating this is used by the DSatur algorithm of Brel\'az~\cite{Brelaz1979} and involves choosing the next node based on the state of the current partial node colouring. Specifically, the \emph{saturation degree} of an uncoloured node $v$ is defined as the number of different colours being used by the neighbours of $v$. The next node to be coloured is then selected using the following rule.
	\begin{definition}[DSatur node selection rule]
		Choose the node with the maximum saturation degree. Break any ties by choosing the node with the largest degree in the subgraph induced by the uncoloured nodes.
	\end{definition}
	This rule seeks to prioritise the colouring of the most ``troublesome'' or ``constrained'' nodes based on the current partial colouring. When used in conjunction with the greedy colouring algorithm, the resultant method runs in $O((m + n) \lg n)$ time and is exact for several graph topologies, including bipartite graphs, trees, cycles, and wheel graphs~\cite{Lewis2021}. Kubale and Jackowski have also found that this rule, when used at Lines~4 and 19 of Algorithm~\ref{alg:bktr}, leads to significantly improved run times compared to other alternatives~\cite{Kubale1985}.
	
	The DSatur algorithm and its backtracking extension are both included in the GCol library. Despite the noted improvements, however, the latter can still exhibit restrictively large run times for non-trivial graphs. For this reason, several metaheuristic-based algorithms are also included in the library. Here, we therefore also make use of a stochastic hybrid evolutionary algorithm (HEA) that was originally proposed by Galinier and Hao~\cite{Galinier1999}. This is among the best-known algorithms for node colouring~\cite{Lewis2021}. The main element of this algorithm is a tabu search process that operates under a fixed number of colours $k$, but that allows clashes to occur (a clash is the occurrence of two adjacent nodes having the same colour). The aim is to alter the colour assignments so that the number of clashes is reduced to zero, at which point $k$ is decremented. The overall evolutionary method maintains a small population of $k$-coloured solutions, which is evolved using selection, recombination, local search, and replacement. In each cycle, two parent solutions are selected from the population, and these are used in conjunction with a specialised recombination operator to produce a new offspring solution. Tabu search is then applied to the offspring for a fixed number of iterations, and the resultant solution is inserted back into the population, replacing its weaker parent. The process continues until a prescribed computation limit or the previously noted lower bound on $\chi(G)$ is achieved.
	
	\begin{figure}[t!]
		\begin{center}
			\includegraphics[trim= 10 10 10 10, clip, scale=0.47]{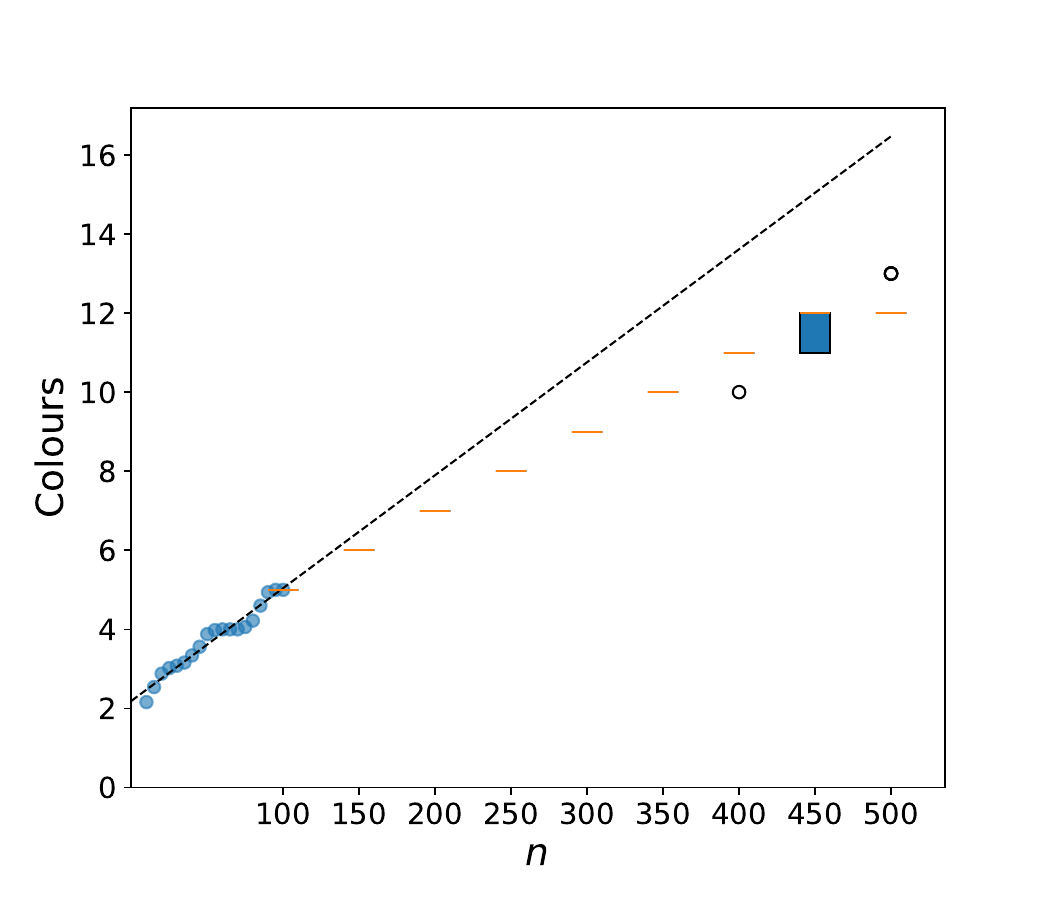}
			\includegraphics[trim= 5 10 10 10, clip, scale=0.47]{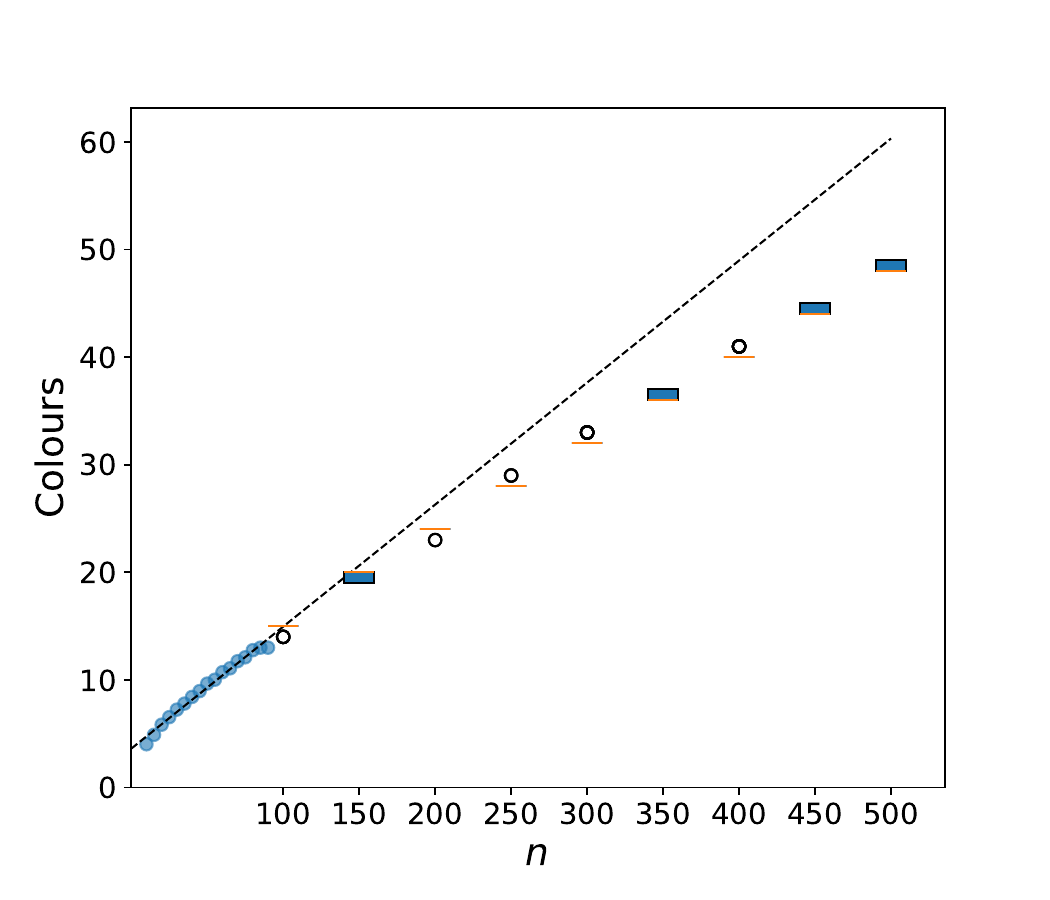}
			\includegraphics[trim= 5 10 10 10, clip, scale=0.47]{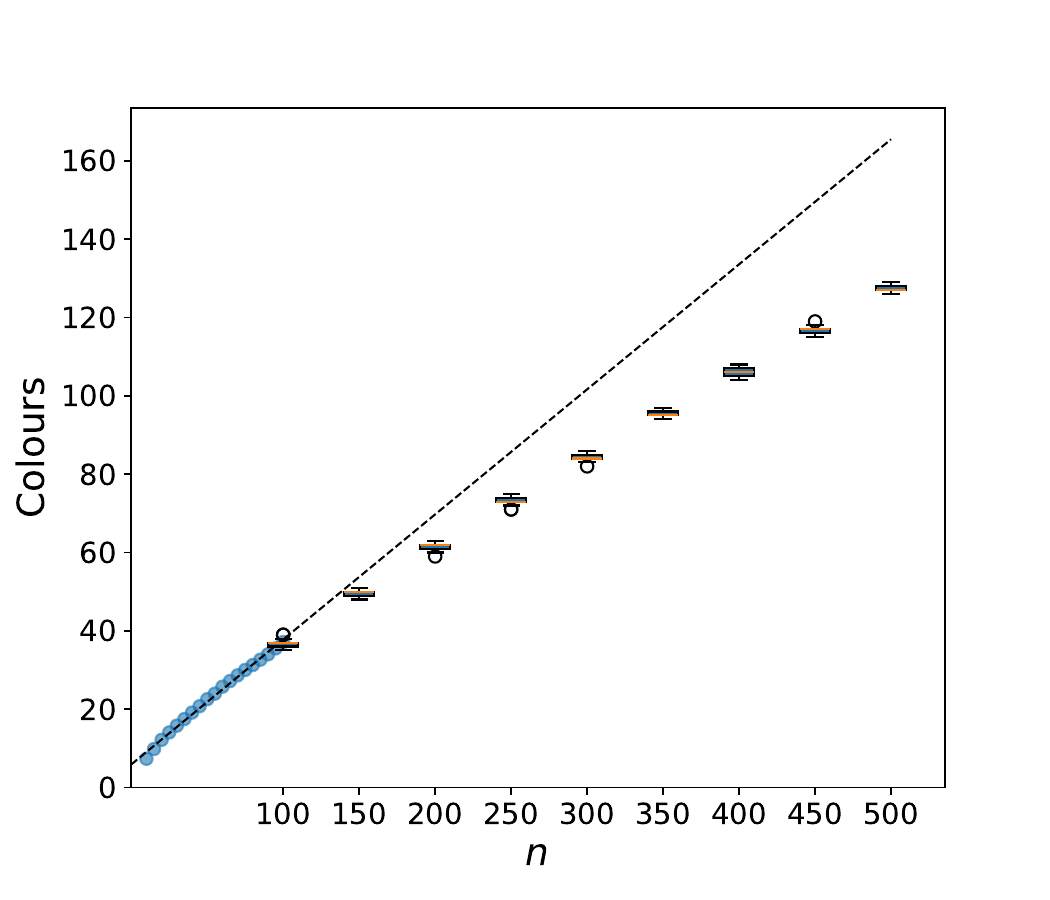}
		\end{center}
		\caption{Performance of the HEA and backtracking algorithms using $G(n,p)$ graphs for $p=0.1$, $0.5$, and $0.9$ respectively. Each point and box plot, generated by the backtracking algorithm and HEA, respectively, was generated from runs on fifty randomly generated graphs.}
		\label{fig:performance}
	\end{figure}
	
	Figure~\ref{fig:performance} demonstrates the performance of the HEA and backtracking algorithms. Here, trials were conducted on sets of randomly generated Erd\H{o}s-Renyi graphs, denoted by $G(n,p)$. Such graphs are constructed by taking a user-prescribed number of nodes $n$ and adding an edge between each node pair at random with probability $p$. The expected number of edges is therefore $\binom{n}{2}p$ and the expected node degree is $(n-1)p$. Here, we show results for differing values for $n$ using $p\in\{0.1,0.5,0.9\}$, representing sparse, medium, and dense graphs, respectively. 
	
	The blue points in Figure~\ref{fig:performance} show the mean chromatic numbers for $G(n,p)$ graphs of up to one hundred nodes. These were found using the backtracking algorithm, which usually completed execution in less than a minute.\footnote{Using a Windows 10 machine with 3.5 GHz quad-core CPUs and 8 GB RAM.} Linear regression lines generated from these points are also shown. The box plots show the number of colours in solutions produced by the HEA, using execution times of up to three minutes. These should only be considered as upper bounds on $\chi(G)$. As indicated, an increase in the number of nodes and/or density tends to lead to larger numbers of colours being needed, though this growth is slightly less than the predicted linear regression. The number of colours needed for each set of $G(n,p)$ graphs also seems to fall within a narrow set of values, indicating stability in algorithm behaviour.

	\subsection{Visualising Node Colourings}
	
	In this section, we illustrate optimal node colourings on a variety of different graphs. When visualising such graphs, observe that the spatial placement of nodes affects interpretability. Good node layouts can reveal structural patterns, clusters, and symmetries, while poor layouts can obscure them. One option in this regard is to use force-directed methods, such as the spring layout~\cite{diBattista1999,Fruchterman1991}. This models nodes as mutually repelling elements, and edges as springs. The method then iteratively adjusts the node positions to minimise an energy function, balancing the attracting forces of edges and the repulsive forces from nodes. The aim is to create an aesthetically pleasing layout where groups of related nodes are close, unrelated nodes are separated, and few edges intersect. 
	
	\begin{figure}[htbp]
		\begin{center}
			\includegraphics[trim= 10 10 10 10, clip, scale=0.44]{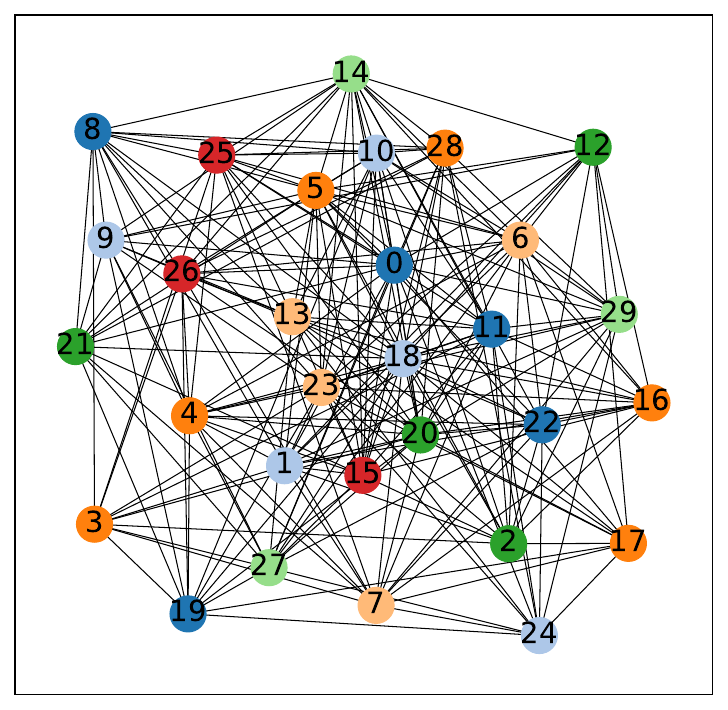}
			\includegraphics[trim= 10 10 10 10, clip, scale=0.44]{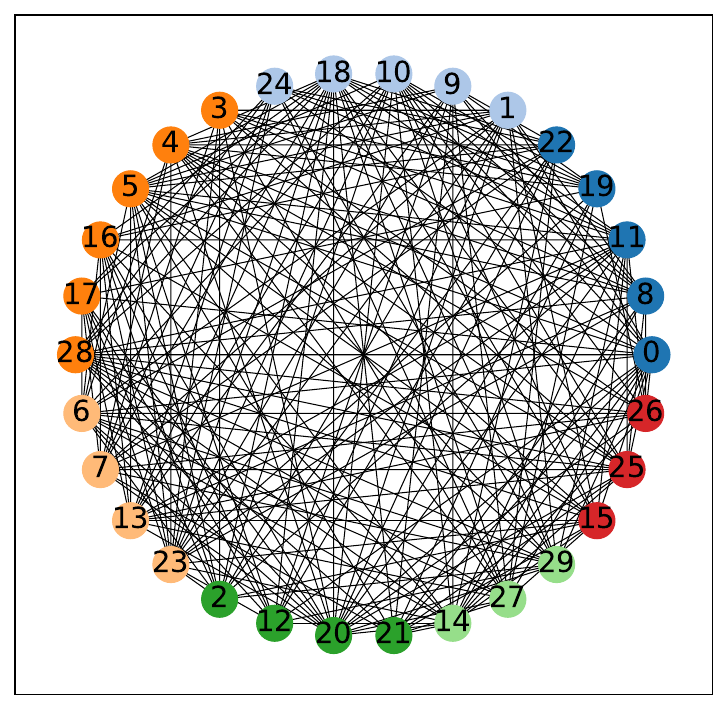}
			\includegraphics[trim= 10 10 10 10, clip,scale=0.44]{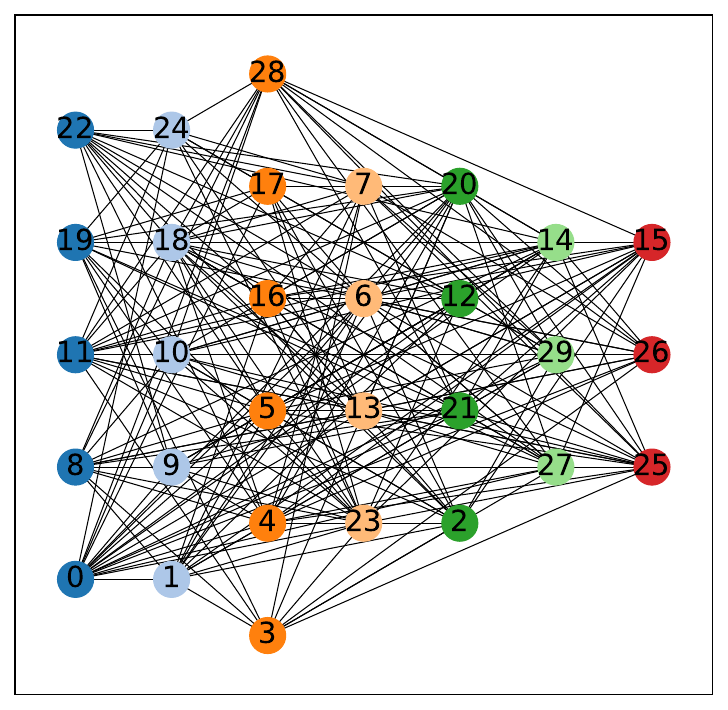}
		\end{center}
		\caption{Three ways of drawing a node colouring of a $G(30,0.5)$ graph. In this case $\chi(G)=7$. Despite looking different, the graph and colouring is identical in each image.}
		\label{fig:random}
	\end{figure}
	
	The leftmost image in Figure~\ref{fig:random} uses a force-directed method to position the nodes of an optimal node colouring of a randomly generated Erd\H{o}s-Renyi graph $G(30,0.5)$. Though this certainly gives a better layout than a random placement, the relatively high number of edges still leads to a rather confusing-looking solution. Two alternative layouts are therefore also shown in the figure. The central image shows all nodes on the circumference of a circle and puts nodes of the same colour in adjacent positions; the second arranges the nodes of each colour into columns. In the latter case, edges can never be vertical, making it easier to see that each column of nodes is an independent set.  
	
	\begin{figure}[htbp]
		\begin{center}
			\includegraphics[trim=10 10 10 10, clip, scale=0.44]{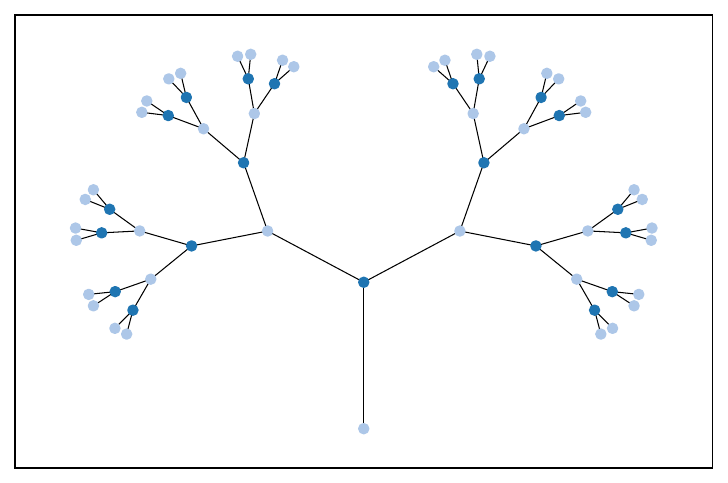}
			\includegraphics[trim=10 10 10 10, clip, scale=0.44]{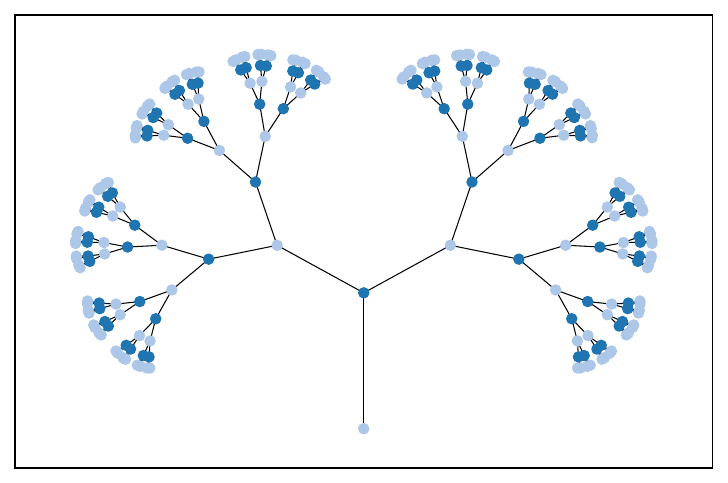}
			\includegraphics[trim=10 10 10 10, clip, scale=0.44]{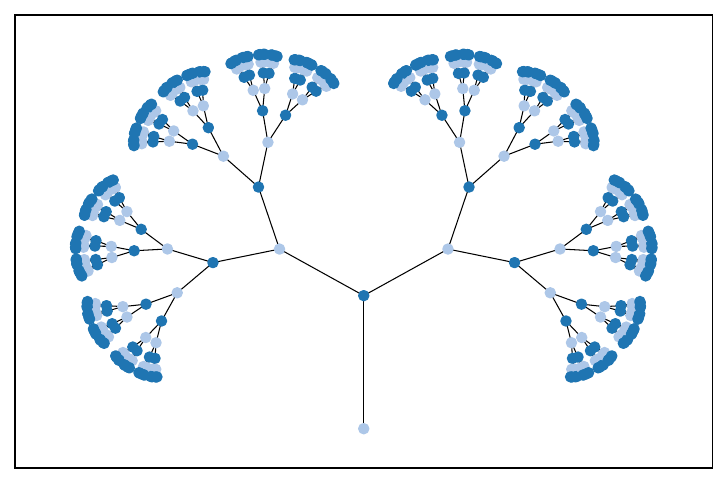}
			\includegraphics[trim=10 10 10 10, clip, scale=0.44]{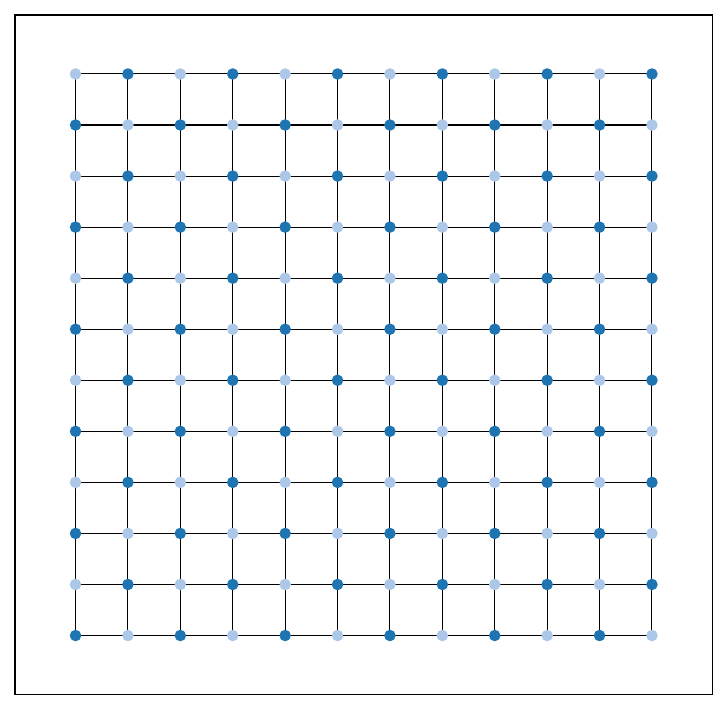}
			\includegraphics[trim=10 10 10 10, clip, scale=0.44]{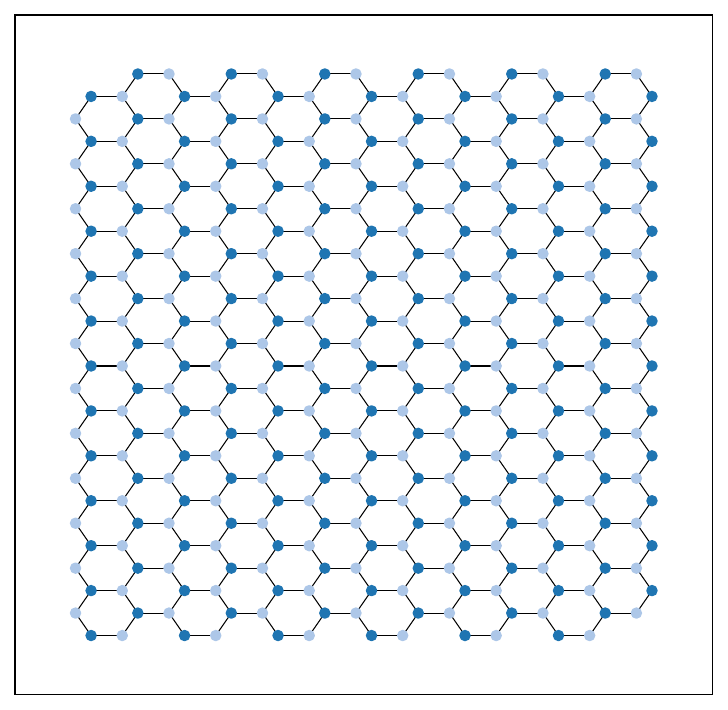}
			\includegraphics[trim=10 10 10 10, clip, scale=0.44]{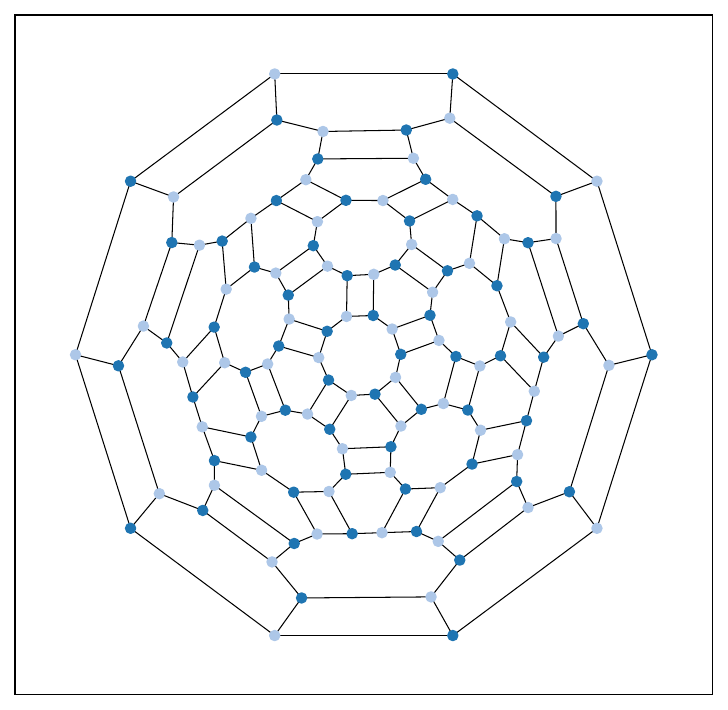}
		\end{center}
		\caption{Optimal node colourings of bipartite graphs. The top three images show binary trees with $64$, $256$, and $512$ nodes; the bottom three show a square lattice, a hexagonal lattice, and the great rhombicosidodecahedral graph (HoG 1122). In all cases $\chi(G)=2$.}
		\label{fig:bipartite-nodes}
	\end{figure}
	
	Node colourings are easier to display when the number of edges and/or colours is small. In some cases, nodes may also have a natural positioning that aids interpretation. Examples of such graphs are shown in Figure~\ref{fig:bipartite-nodes}. The top three images show fractal layouts of binary trees and are somewhat reminiscent of Chinese willow patterns. The bottom three show embeddings of a square lattice, a hexagonal lattice, and the great rhombicosidodecahedral graph. Note that none of these graphs contains an odd-length cycle. As such, they are bipartite and feature $\chi(G)=2$. This means that they are also solved by the DSatur algorithm without the need for any additional optimisation.  
	
	\begin{figure}[htbp]
		\begin{center}
			\includegraphics[trim=10 10 10 10, clip, scale=0.44]{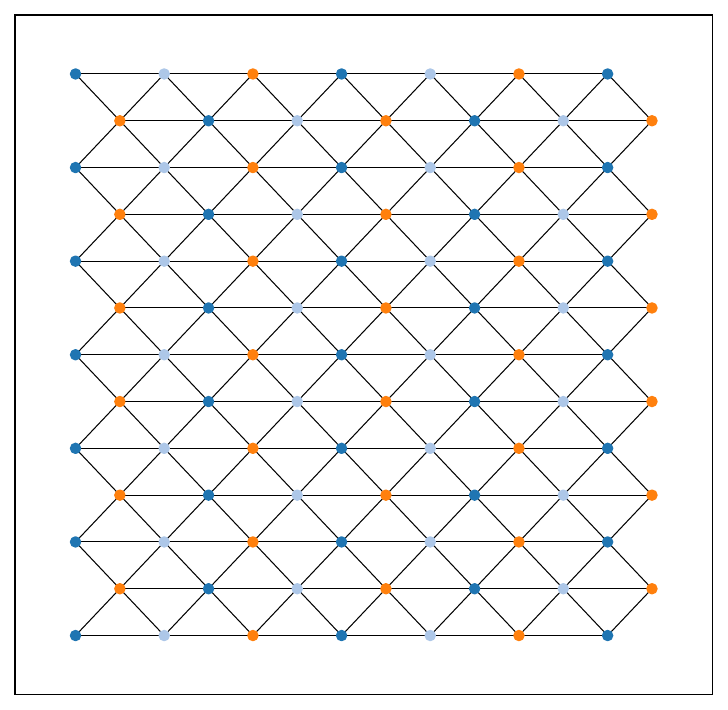}
			\includegraphics[trim=10 10 10 10, clip, scale=0.44]{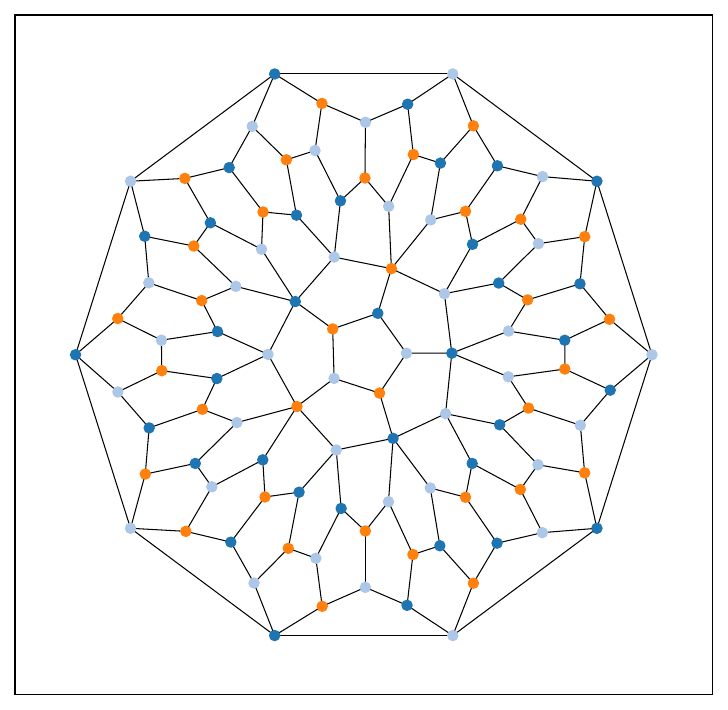}
			\includegraphics[trim=10 10 10 10, clip, scale=0.44]{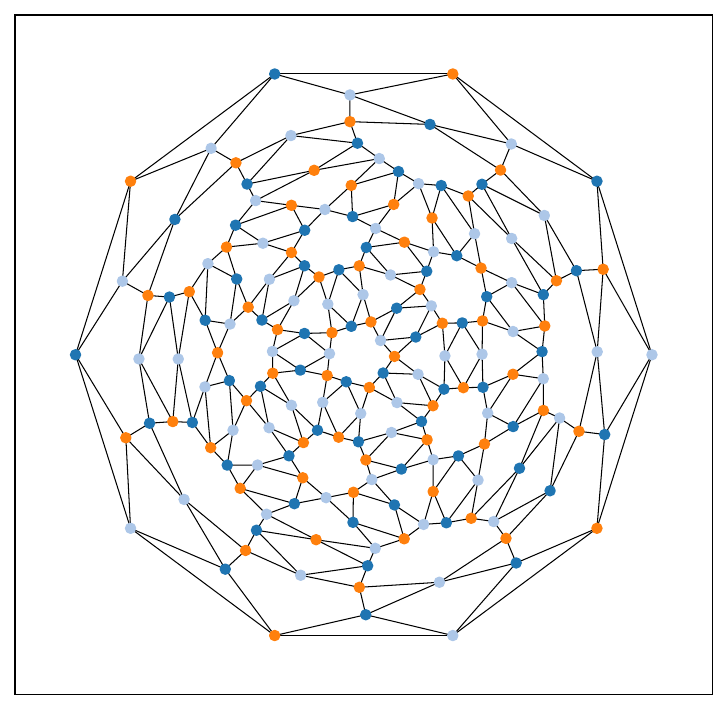}
		\end{center}
		\caption{Optimal node colourings of a triangular lattice graph, Thomassen graph (HoG 1347), and the great rhombicosidodecahedral line graph (HoG 51392). In all cases $\chi(G)=3$.}
		\label{fig:optimal-nodes}
	\end{figure}
	
	Figure~\ref{fig:optimal-nodes} shows optimal node colourings for a triangular lattice, a Thomassen graph on $105$ nodes, and the great rhombicosidodecahedral line graph. By visual inspection, we see that these graphs contain odd-length cycles, implying that they are not bipartite and, therefore, $\chi(G)\geq 3$. In fact, $\chi(G)=3$ for each of these graphs, as demonstrated in the figure. Here, optimal solutions for the first two examples were again found using the DSatur algorithm, whereas the third graph required a small amount of additional optimisation.
	
	\section{Edge Colouring}
	\label{sec:edgecol}
	
	Recall from Definition~\ref{def:edgecol} that an edge colouring is an assignment of colours to edges so that adjacent edges have different colours, with the aim of minimising the number of colours. The set of edges assigned to a particular colour corresponds to a matching; hence, an equivalent objective is to partition the graph's edges into a minimum number of matchings.
	
	In edge colourings, the edges that end at any particular node must always have different colours. Hence, for any graph $G$, $\Delta(G)\leq\chi'(G)$, where $\Delta(G)$ denotes the maximum degree in $G$. For bipartite graphs, K\"onig's theorem tells us that $\Delta(G)$ is always sufficient. Vizing's theorem generalises this result, stating that, for any graph $G$, the chromatic index is either $\Delta(G)$ or $\Delta(G)+1$. This means that all graphs can be classified into two classes: Class~1, in which $\chi'(G)=\Delta(G)$; and Class 2, where $\chi'(G)=\Delta(G)+1$. In general, determining if a graph belongs to Class~1 is $\mathcal{NP}$-complete~\cite{Hoyler1981}, though some topologies like bipartite and complete graphs can be solved in polynomial time~\cite{Kirkman1847}. Edge colourings using $\Delta(G)+1$ colours can also be determined in polynomial time for any graph topology~\cite{Misra1992}.
	
	\begin{figure}[htbp]
		\begin{center}
			\includegraphics[trim=10 10 10 10, clip, scale=0.44]{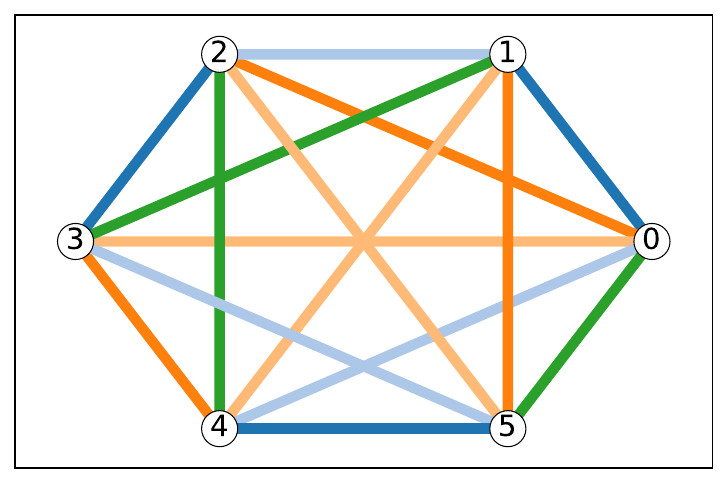}
			\includegraphics[trim=10 10 10 10, clip, scale=0.44]{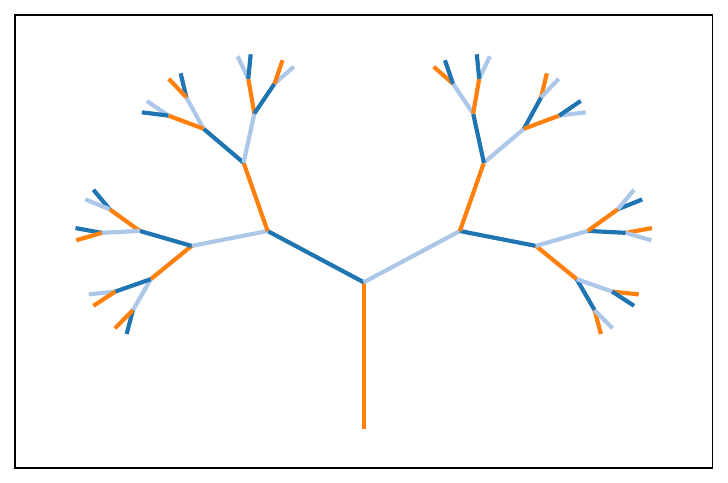}
			\includegraphics[trim=10 10 10 10, clip, scale=0.44]{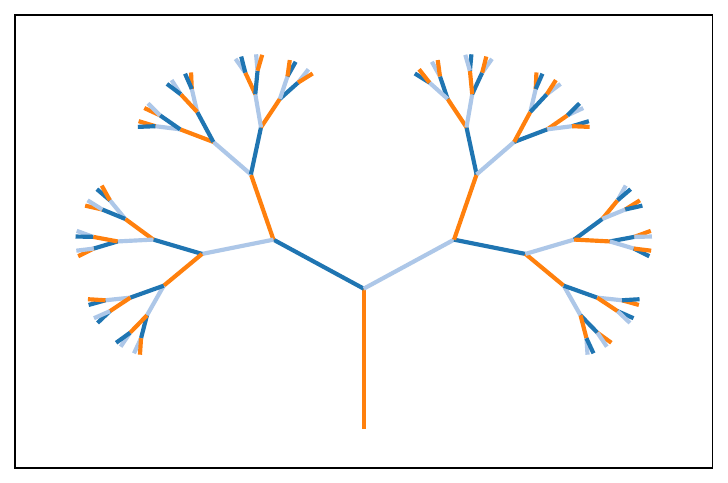}
			\includegraphics[trim=10 10 10 10, clip, scale=0.44]{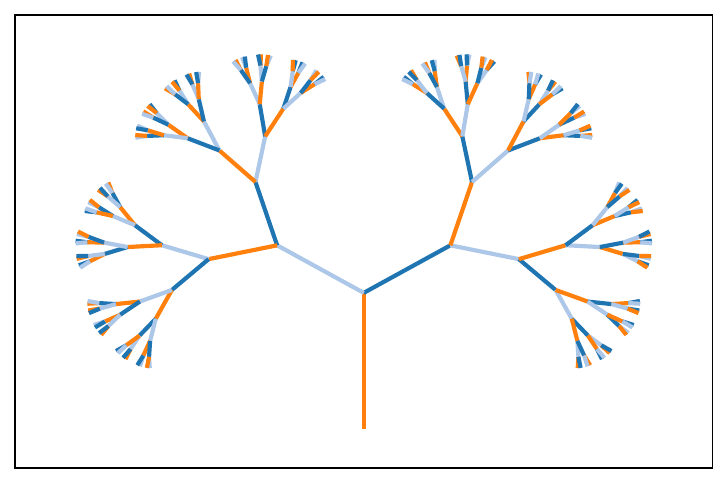}
			\includegraphics[trim=10 10 10 10, clip, scale=0.44]{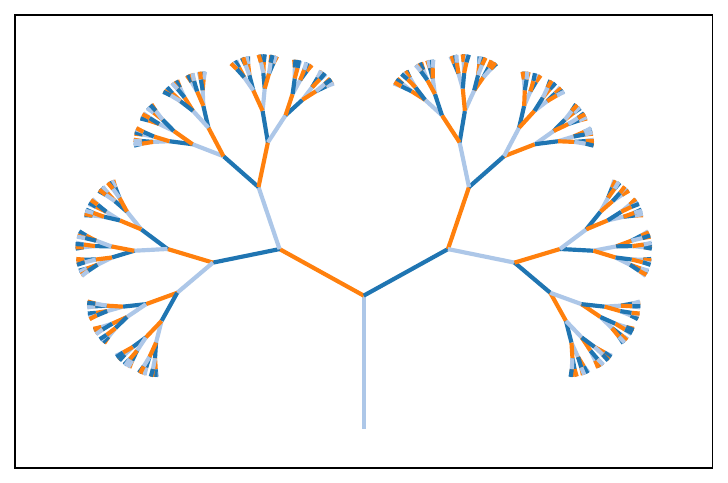}
			\includegraphics[trim=10 10 10 10, clip, scale=0.44]{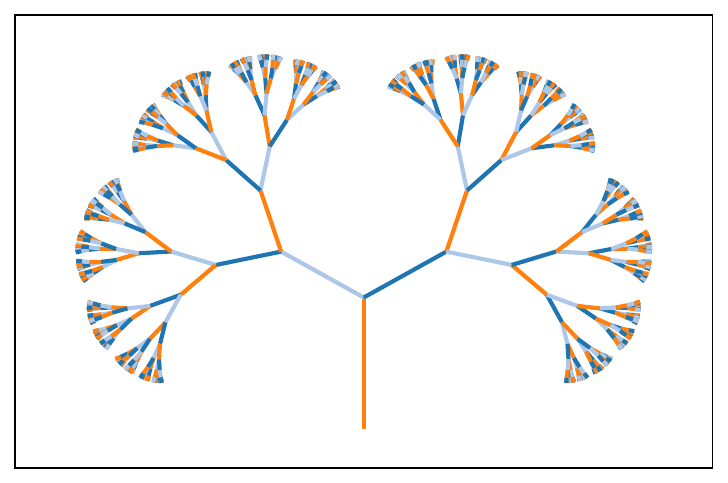}
		\end{center}
		\caption{Optimal edge colourings for a complete graph (top left), and binary trees with 64, 128, 256, 512, and 1024 nodes, respectively.}
		\label{fig:edgcoltrees}
	\end{figure}
	
	Edge colouring has applications in the construction of sports leagues, where a set of teams are required to play each other over a series of rounds~\cite{deWerra1988}. The first image in Figure~\ref{fig:edgcoltrees} shows a complete graph on six nodes, one node per team. Here, edges represent matches between teams, and each colour gives a single round in the schedule. Hence, the ``light blue'' round involves matches between Teams 1 and 2, 0 and 3, and 4 and 5, for example. Complete graphs are known to belong to Class~1 if and only if $n$ is even.
	
	Figure~\ref{fig:edgcoltrees} also shows edge colourings for the binary tree fractals seen earlier. In these images, we do not display graph nodes to keep visual clutter to a minimum. Following K\"onig's theorem, $\Delta(G)=3$ colours are sufficient in these cases. It is interesting to note in these images that we can describe the (unique) path between the root note and any other node by simply listing the observed sequence of colours. In fact, a generalisation of this property exists for all graphs: in an edge colouring of a graph $G$, every walk starting at a node $v$ is describable by a tuple $(v, C)$, where $C=(c_1, c_2, \ldots, c_l)$ is a sequence in which each $c_i$ gives the colour of the $i$th edge in the corresponding walk. 
	
	\begin{figure}[htbp]
		\begin{center}
			\includegraphics[trim=30 30 30 30, clip, scale=0.68]{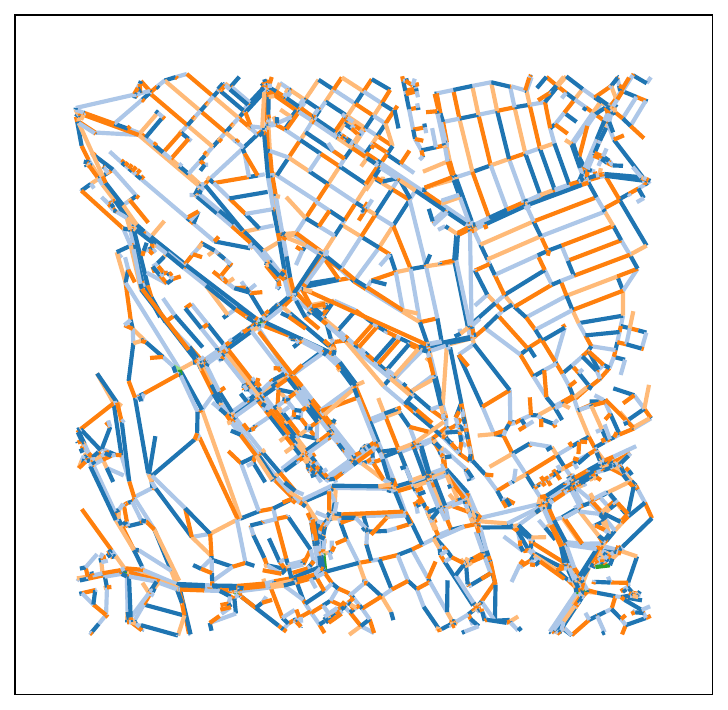} \hspace{0.75cm}
			\includegraphics[trim=30 30 30 30, clip, scale=0.68]{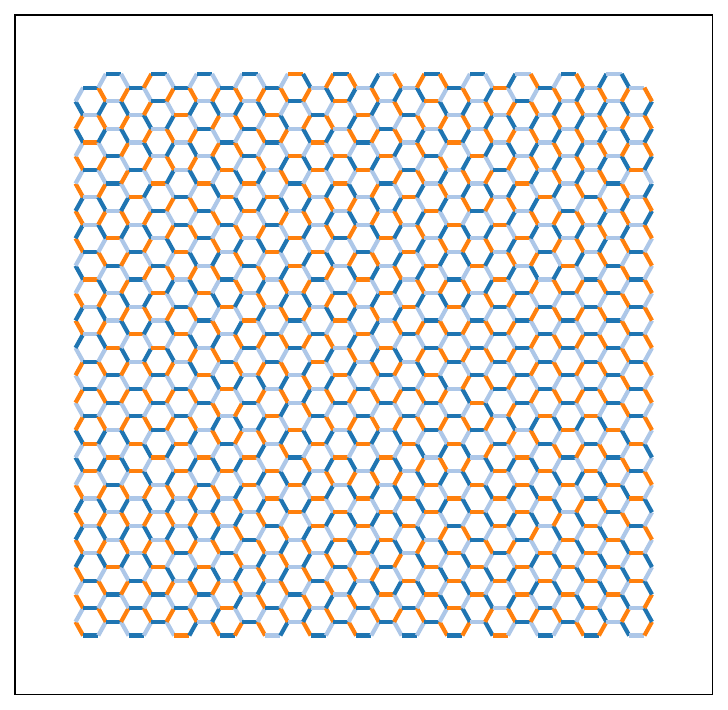}
		\end{center}
		\caption{An optimal edge colouring of the street map of central Cardiff, Wales (left) and the hexagonal lattice graph (right). The chromatic indices of these graphs are six and three, respectively. In both cases, $\chi'(G)=\Delta(G)$.}
		\label{fig:edgcolmaps}
	\end{figure}
	
	This property is exemplified by the street map shown in Figure~\ref{fig:edgcolmaps}, where $\chi'(G)=\Delta(G)=6$. In GPS navigation systems, routes between two locations are defined using the start node and a sequence of tags that identify the street segments followed. OpenStreetMap, for example, uses a 64-bit integer tag for each street segment in the world map. On the other hand, just six colours are needed with this map, meaning that routes can be defined using the starting node and just three bits per edge. 
	
	\begin{figure}[htbp]
		\begin{center}
			\includegraphics[trim=10 10 10 10, clip, scale=0.44]{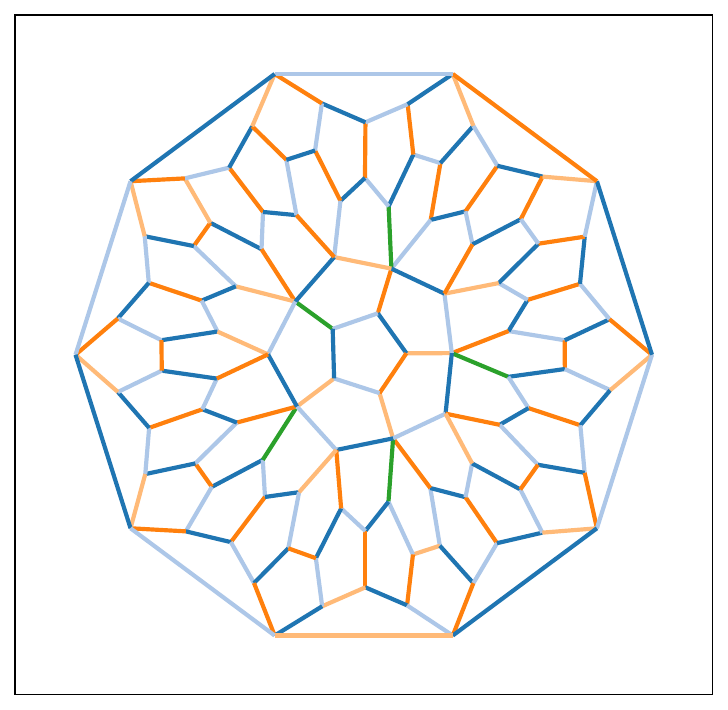}
			\includegraphics[trim=10 10 10 10, clip, scale=0.44]{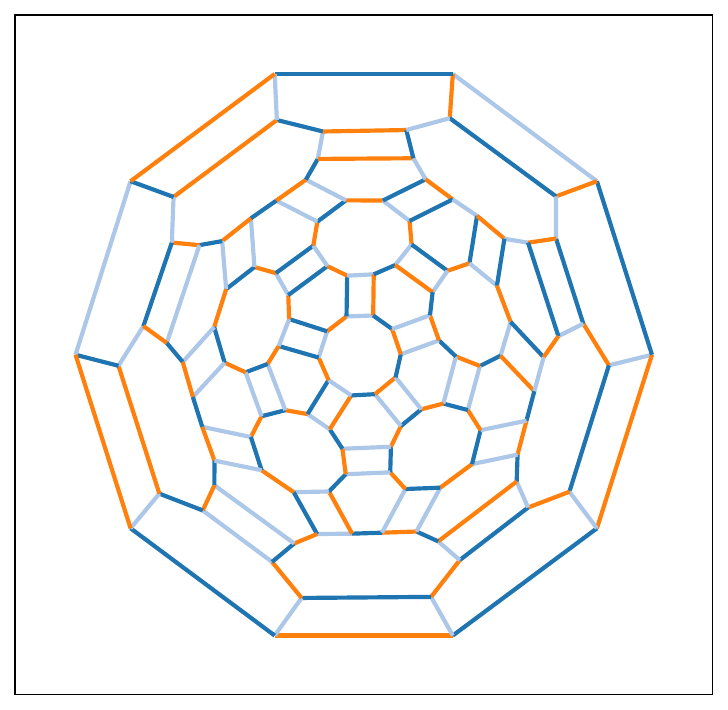}
			\includegraphics[trim=10 10 10 10, clip, scale=0.44]{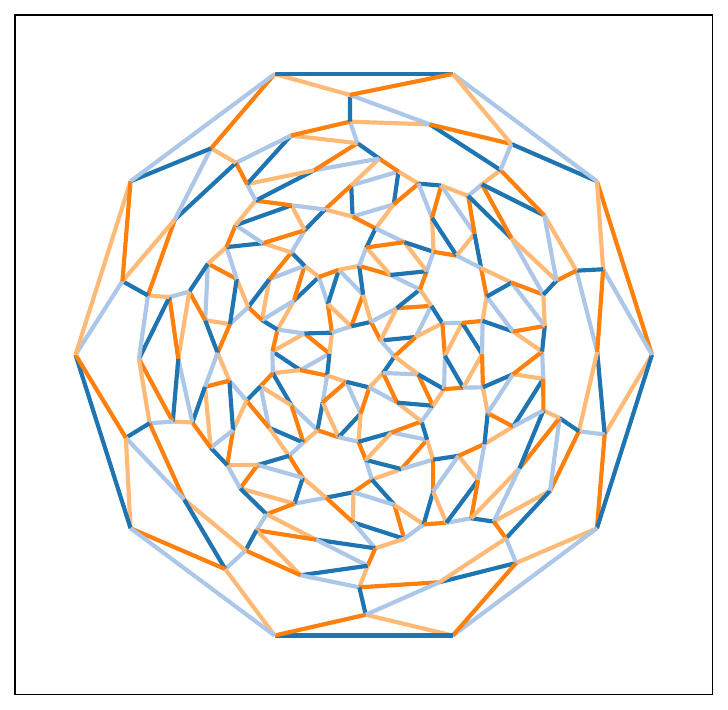}
		\end{center}
		\caption{Optimal edge colourings of the Thomassen graph, the great rhombicosidodecahedral graph, and its line graph. The chromatic indices of these graphs are five, three, and four, respectively. In all cases, $\chi'(G)=\Delta(G)$.}
		\label{fig:moreedgecols}
	\end{figure}
	
	Finally, Figure~\ref{fig:moreedgecols} shows optimal edge colourings of three further Class~1 graphs. These examples are reminiscent of crochet doily patterns or, perhaps, Ojibwe dream catchers.
	
	\section{Face Colouring}
	\label{sec:facecol}
	
	The fact that all planar embeddings feature a face chromatic number $\chi_f(G)\leq 4$ was first conjectured in 1852 by Francis Guthrie. Surprisingly, however, it would take over a hundred years---and the research endeavours of many notable mathematicians including Augustus De~Morgan, William Hamilton, Arthur Cayley and Percy Heawood---before this simple sounding proposition was eventually proved by Appel and Haken in 1977~\cite{Appel1977,Wilson2003}. 
	
	\begin{figure}[htbp]
		\begin{center}
			\includegraphics[trim=10 10 10 10, clip, scale=0.44]{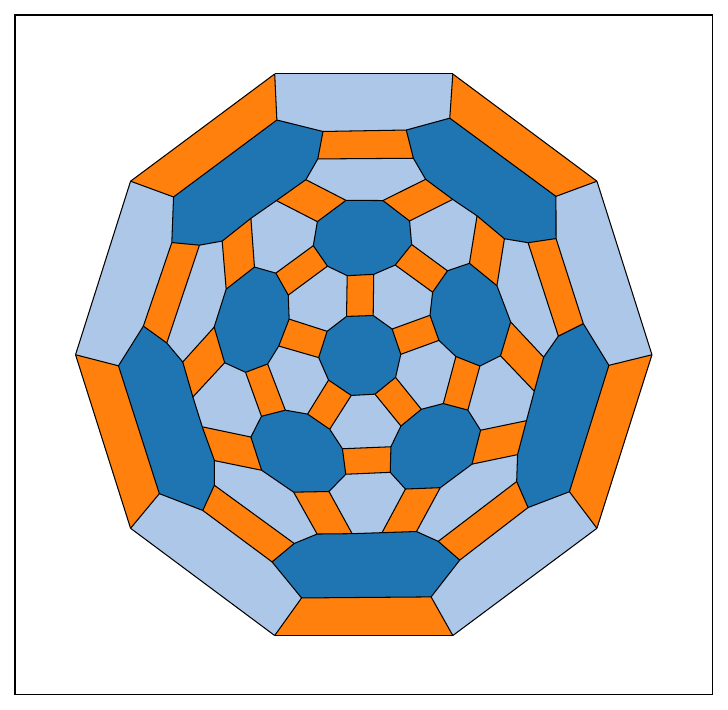}
			\includegraphics[trim=10 10 10 10, clip, scale=0.44]{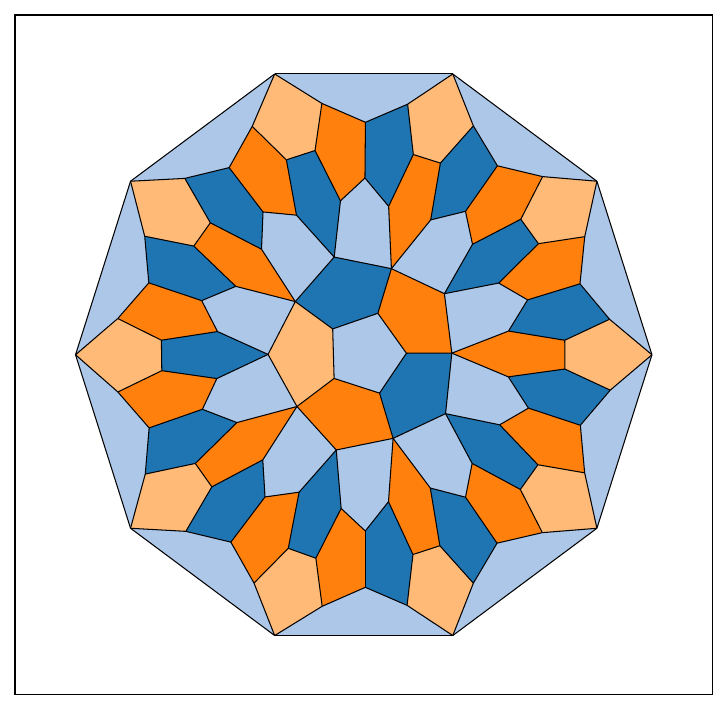}
			\includegraphics[trim=10 10 10 10, clip, scale=0.44]{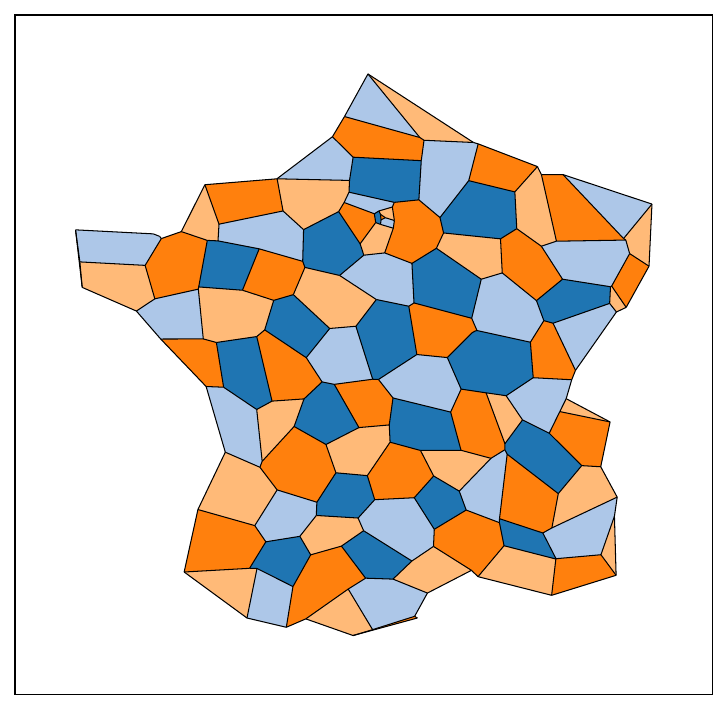}
		\end{center}
		\caption{Optimal face colourings of the great rhombicosidodecahedral graph, the Thomassen graph, and a map of the administrative departments of France. Here, the face chromatic numbers are three, four, and four, respectively.}
		\label{fig:facefourcol}
	\end{figure}
	
	Figure~\ref{fig:facefourcol} shows three example face colourings. For ease of viewing, nodes are again not shown here, but their presence should be assumed whenever edges are seen to meet. Note that in each of these images, the unbounded face can be coloured dark blue, maintaining the number of colours. The first two solutions in this figure were found using the DSatur algorithm, while the third also required a small amount of optimisation. In this figure, the central face of the Thomassen graph illustrates why four colours are often needed. As shown, this central face is adjacent to five surrounding faces. Together, these five faces form an odd cycle, necessarily requiring three different colours, so the central face must then be allocated to a fourth colour. 
	
	\begin{figure}[b!]
		\begin{center}
			\includegraphics[trim=10 10 10 10, clip, scale=0.44]{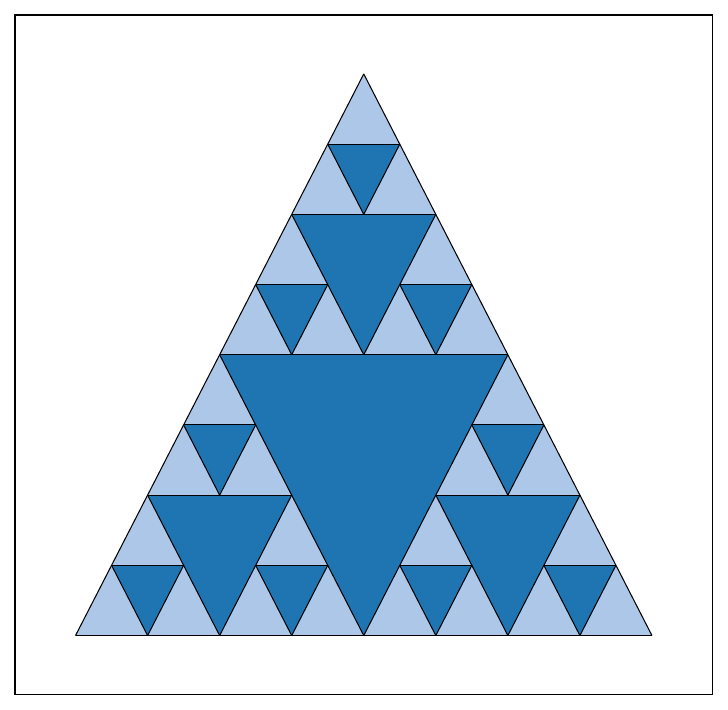}
			\includegraphics[trim=10 10 10 10, clip, scale=0.44]{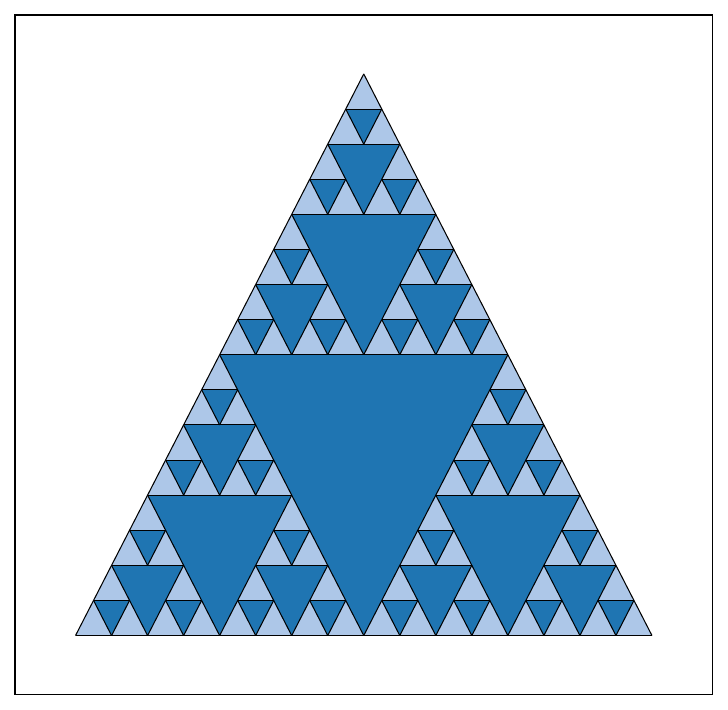}
			\includegraphics[trim=10 10 10 10, clip, scale=0.44]{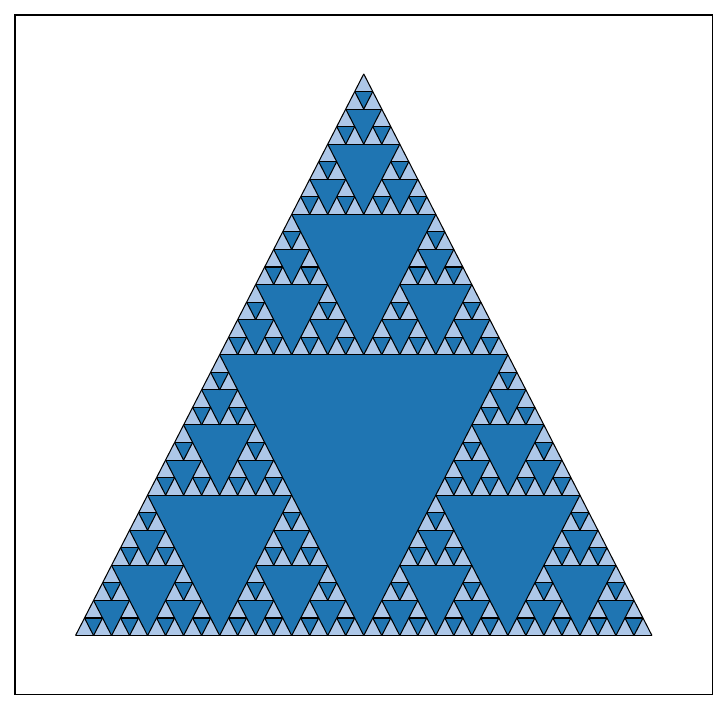}
			\includegraphics[trim=10 10 10 10, clip, scale=0.44]{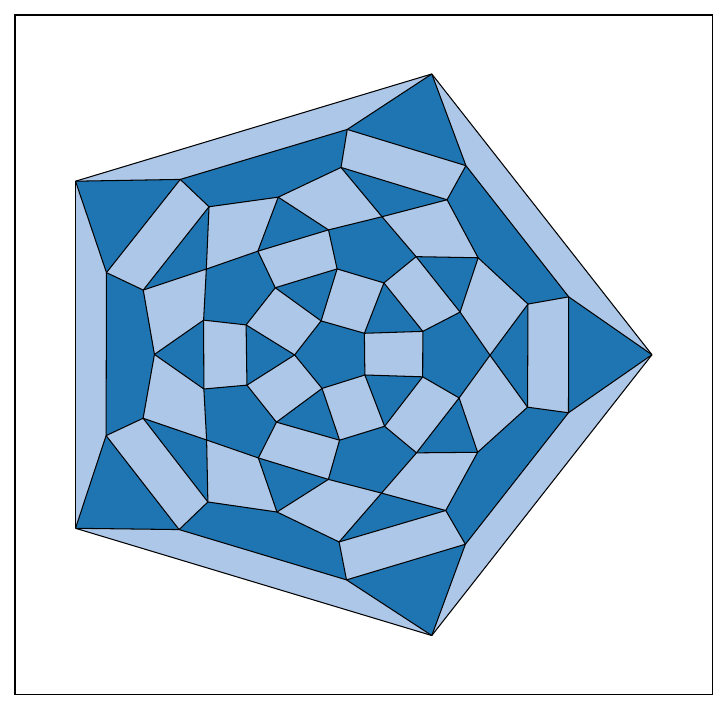}
			\includegraphics[trim=10 10 10 10, clip, scale=0.44]{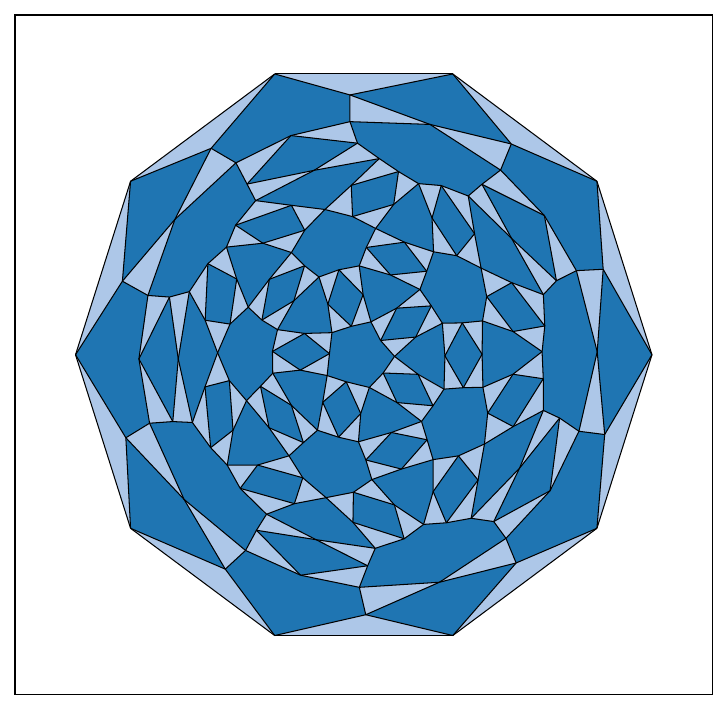} 
			\includegraphics[trim=10 10 10 10, clip, scale=0.44]{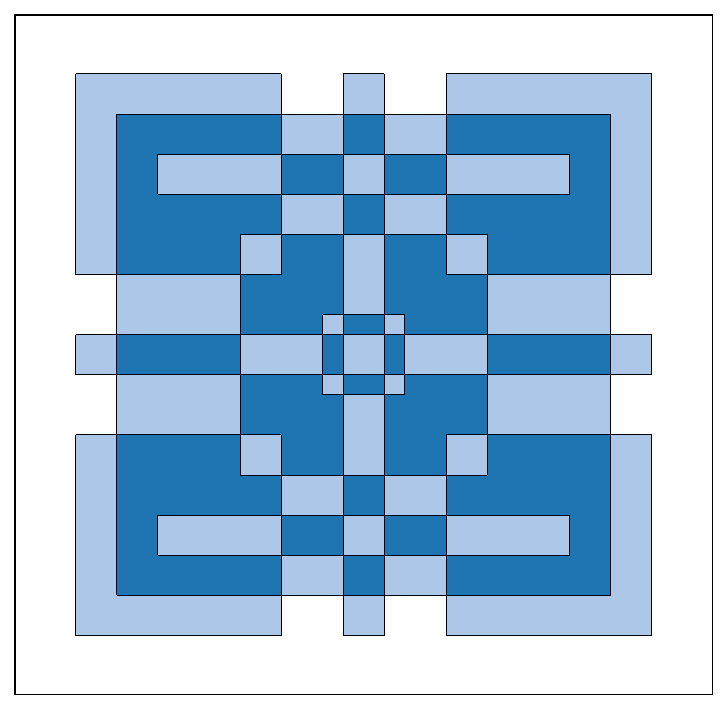} 
		\end{center}
		\caption{Illustrations of how Eulerian planar graphs always feature $\chi_f(G)=2$. The first three examples show the Sierpinski triangle fractal at 3, 4, and 5 levels of recursion; the next two examples show the small rhombicosidodecahedral graph (HoG 1317) and the great rhombicosidodecahedral line graph. The final example is formed by overlaying an arbitrary set of closed curves (rectangles in this case). Here, whenever the curves intersect, we form a node of even degree, making the graph Eulerian.}
		\label{fig:eulerfacecol}
	\end{figure}
	
	In graph theory, a graph is called \emph{Eulerian} if and only if it is connected and the degrees of all its nodes are even. A well-known theorem in this regard is the following.
	\begin{theorem}
		\label{thm:planar}
		A connected planar graph $G$ is Eulerian if and only if its dual graph $G^*$ is bipartite.
	\end{theorem}
	Because a face colouring of a planar graph $G$ corresponds to a node colouring of $G^*$, this theorem implies that Eulerian planar graphs always feature a face chromatic number $\chi_f(G)=\chi(G^*)=2$. Examples of this are shown in Figure~\ref{fig:eulerfacecol} where, as required, all nodes have an even degree. In addition, the unbounded faces might also be coloured dark blue as before. Practical examples of Theorem~\ref{thm:planar} can be seen in chess boards, spirograph patterns, and many forms of Islamic art, all of which feature underlying graphs that are both planar and Eulerian. Common tiling patterns involving square, rectangular, or triangular tiles are also characterised by such graphs, as seen in the well-known ``chequered'' tiling style.
	
	\begin{figure}[htbp]
		\begin{center}
			\includegraphics[trim=10 10 10 10, clip, scale=0.68]{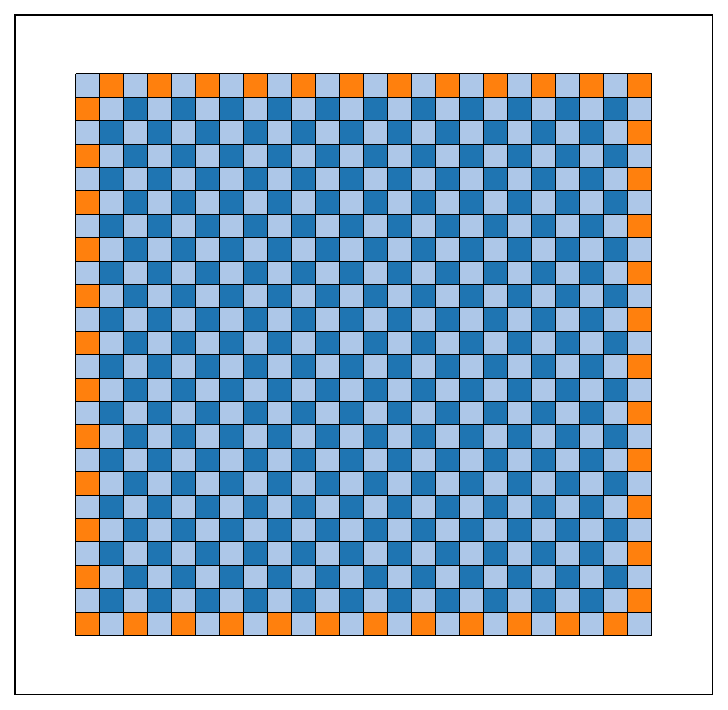}
			\includegraphics[trim=10 10 10 10, clip, scale=0.68]{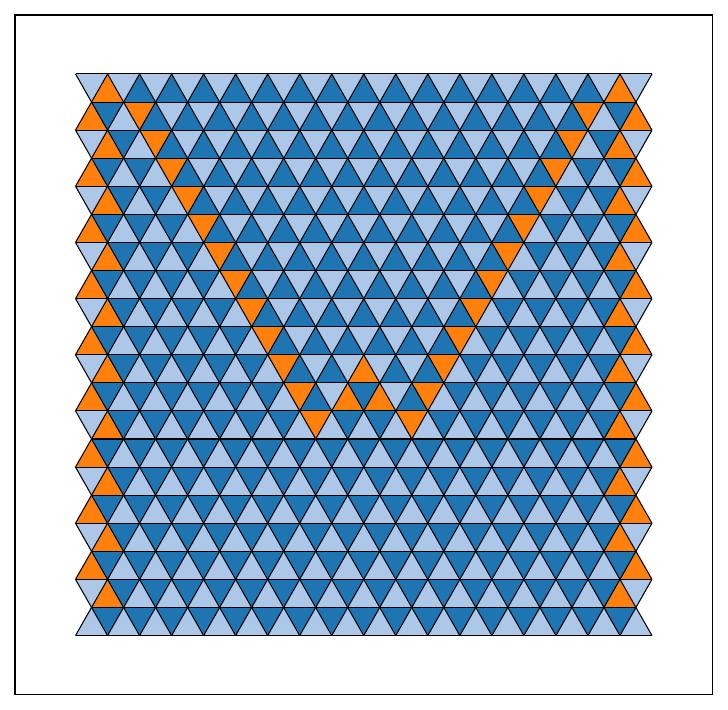}
			\includegraphics[trim=10 10 10 10, clip, scale=0.68]{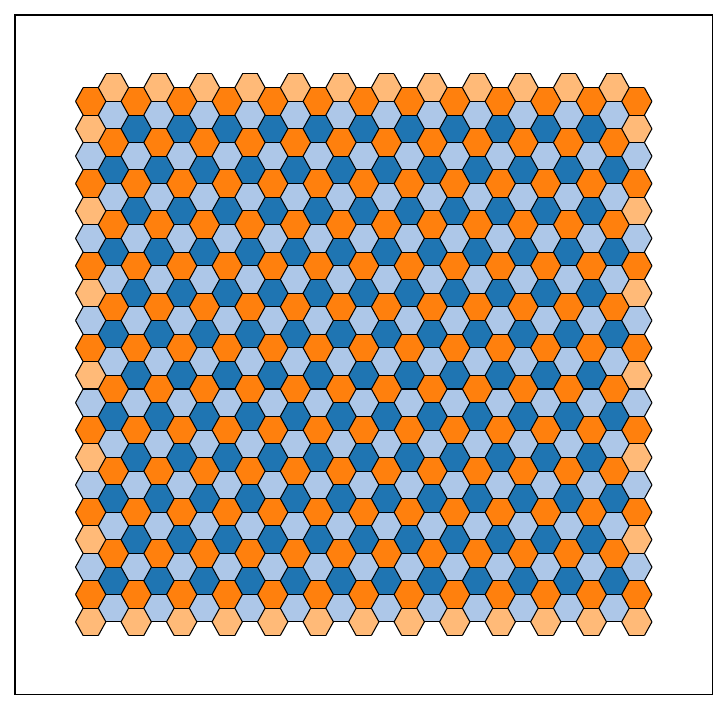}
			\includegraphics[trim=10 10 10 10, clip, scale=0.68]{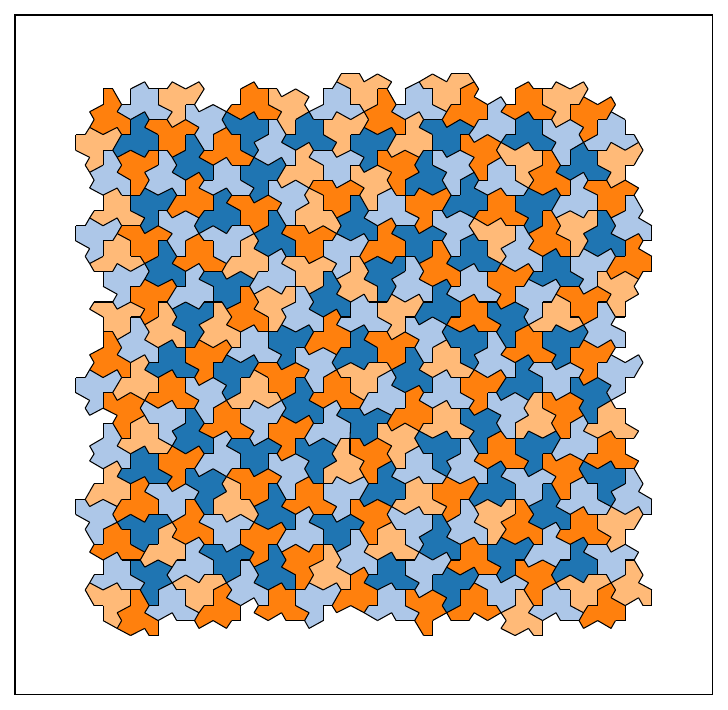}
		\end{center}
		\caption{The first three images show optimal face colourings of the square lattice, triangular lattice, and hexagonal lattice, including the unbounded face. The fourth example shows an optimal colouring of the aperiodic tiling pattern formed by the ``hat'' tile~\cite{Smith2024}.}
		\label{fig:tiling}
	\end{figure}
	
	Four example tiling patterns are shown in Figure~\ref{fig:tiling}. The first two examples show formations of square and triangular tiles and most faces have assumed one of two colours in a periodic formation; however, note that some of the nodes in the graphs' boundaries are odd in degree, meaning that a third colour is needed to allow the unbounded face to also be assigned to dark blue.\footnote{When tiling a wall or floor, we tend not to worry about the notion of an unbounded face, however.} The third example shows a similar result for hexagonal tile shapes. This uses a repeating pattern of three colours in the main body, with an additional colour needed at the boundary to deal with the unbounded face. In contrast, the fourth example shows a recently-discovered \emph{aperiodic} tiling pattern~\cite{Smith2024}. This requires four colours, which are allocated in a much more haphazard fashion. 
	
	\begin{figure}[htbp]
		\begin{center}
			\includegraphics[trim=57 50 50 57, clip, scale=0.645]{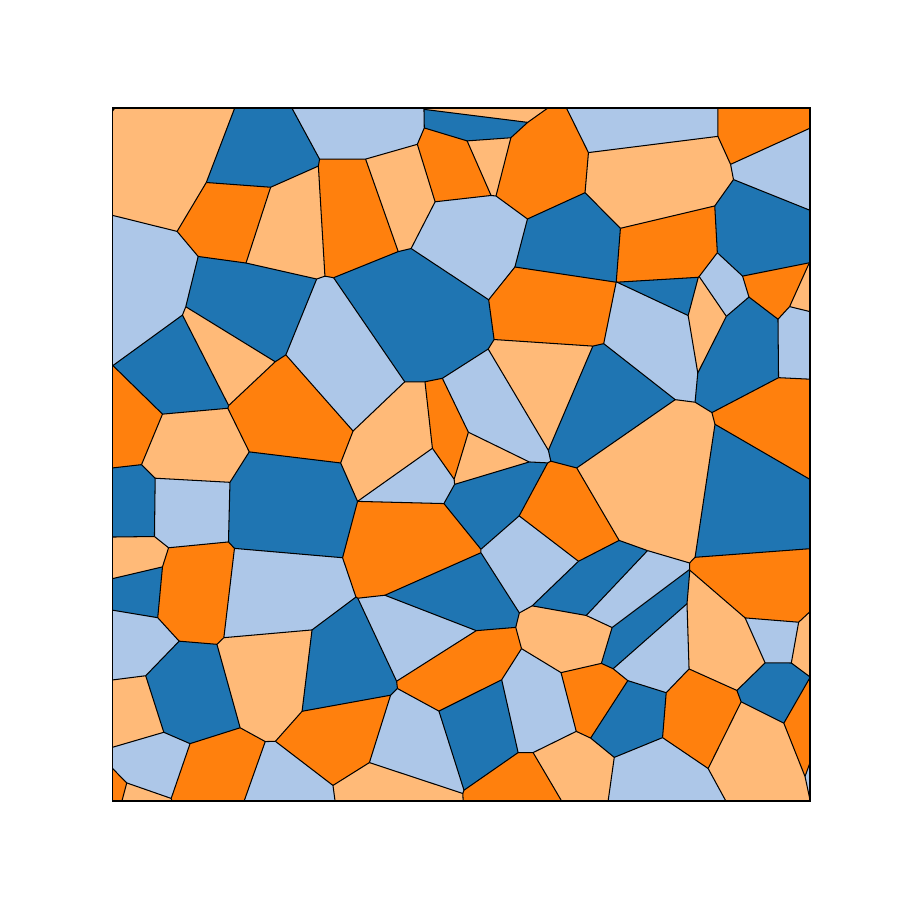} 
			\hspace{0.75cm}
			\includegraphics[trim=37 37 37 37, clip, scale=0.787]{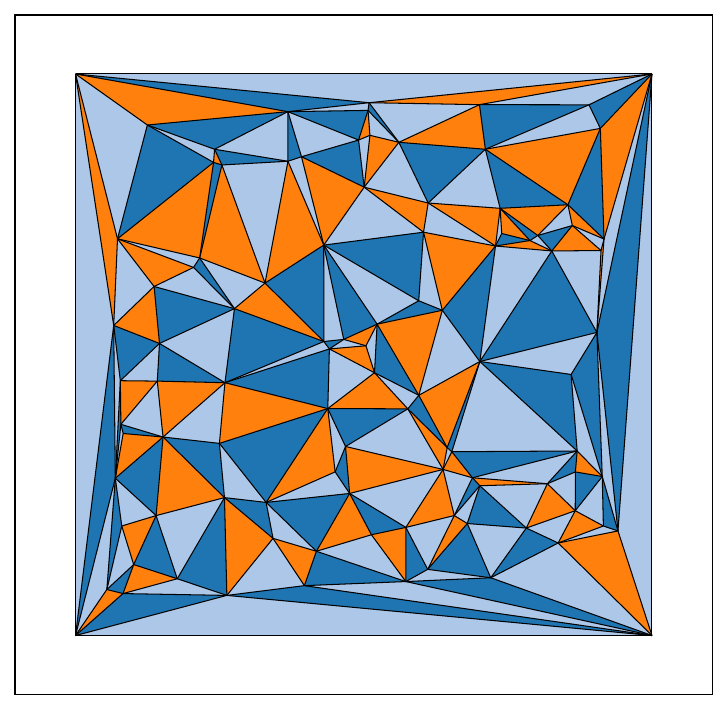}
		\end{center}
		\caption{The left image shows an optimal face colouring of a Voronoi diagram. This was formed from 100 seed points randomly placed in the unit square. The right image shows an optimal face colouring of its dual--- that is, a Delaunay triangulation formed from the same set of points. The face chromatic numbers are four and three, respectively.}
		\label{fig:voronoi}
	\end{figure}
	
	Our next example considers Voronoi diagrams. Given a prescribed set of ``seed'' points $P$ on the plane, a Voronoi diagram is a partition of the plane into ``cells'' so that each cell contains exactly one member $p\in P$, and all other points in the cell are closer to $p$ than any other member of $P$. The cells of a Voronoi diagram cannot overlap; as such, their boundaries collectively define a planar embedding whose faces can be coloured. An example is shown in Figure~\ref{fig:voronoi}.
		
	A set of seed points $P$ can also be used to generate a Delaunay triangulation. These are a triangulation of the convex hull of $P$ such that the circumcircle of every triangle contains no other points of $P$ in its interior. Equivalently, two points $p,q\in P$ are joined in a Delaunay triangulation if and only if there exists a circle passing through $p$ and $q$ that contains no other points of $P$. The Delaunay triangulation of a set $P$ is the dual of the Voronoi diagram of $P$ and, as such, is also a planar graph. Figure~\ref{fig:voronoi} shows an example face colouring of a Delaunay triangulation. Note that this graph has a face chromatic number of three, which is a common feature of Delaunay triangulations, though this is by no means guaranteed.
	
	\begin{figure}[htbp]
		\begin{center}
			\includegraphics[trim=0 0 0 0, clip, scale=0.21]{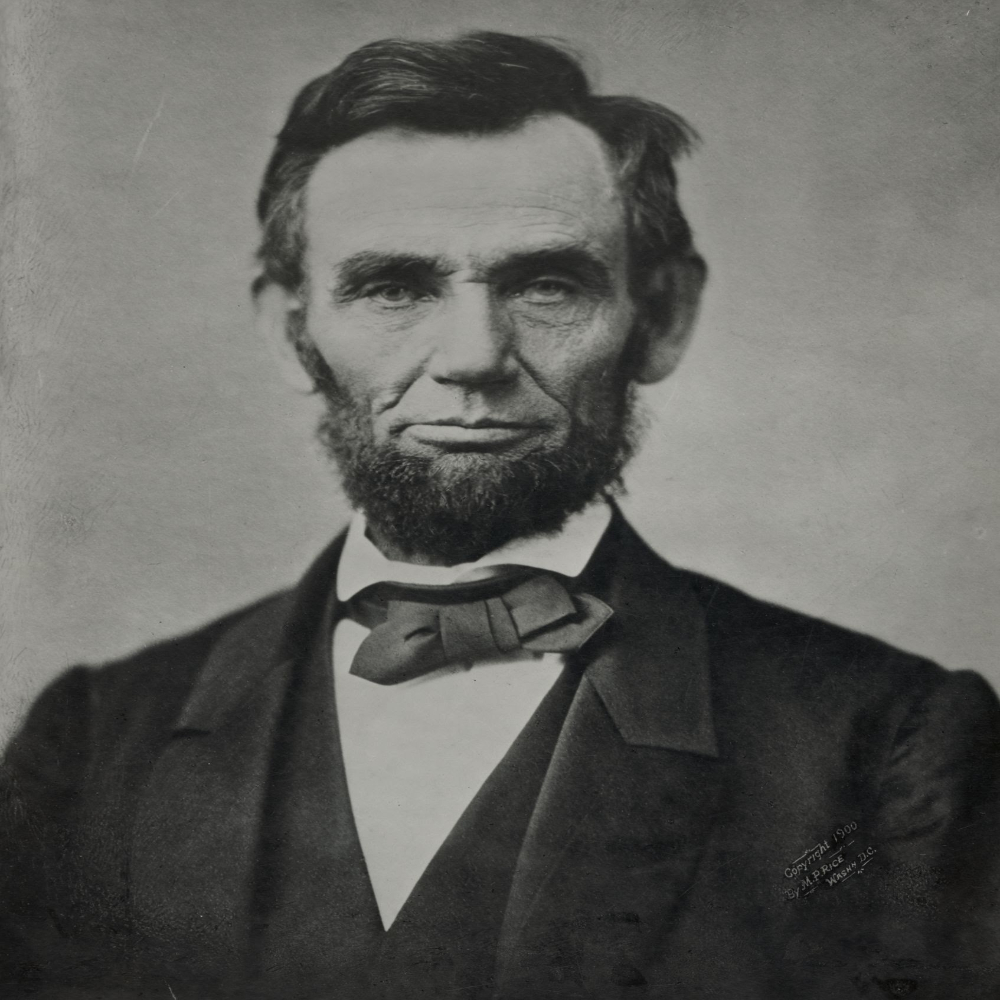} 
			\hspace{0.75cm}
			\includegraphics[trim=37 37 37 37, clip, scale=0.787]{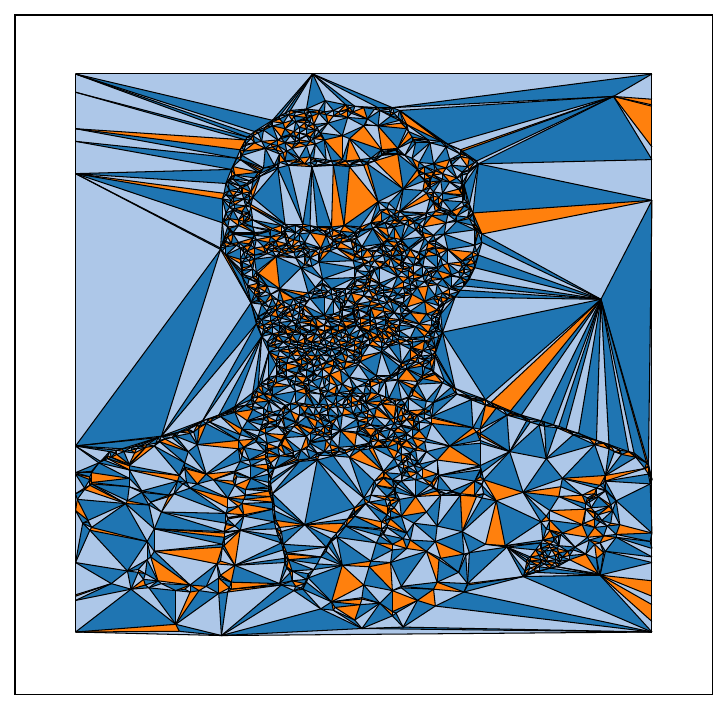}
		\end{center}
		\caption{A Delaunay triangulation (and optimal face colouring) generated from the 1863 portrait of Abraham Lincoln by Alexander Gardner.}
		\label{fig:lincoln}
	\end{figure}
		
	We can, of course, form Delaunay triangulations from any set of points. In Figure~\ref{fig:lincoln}, for example, we have used the online tool of Fischer~\cite{Fischer2026} to show how an arbitrary image can be triangulated and converted into planar graphs. We could use such a process to produce ``guess the identity'' puzzles. 
	
	\begin{figure}[htbp]
		\begin{center}
			\includegraphics[trim=10 10 10 10, clip, scale=0.68]{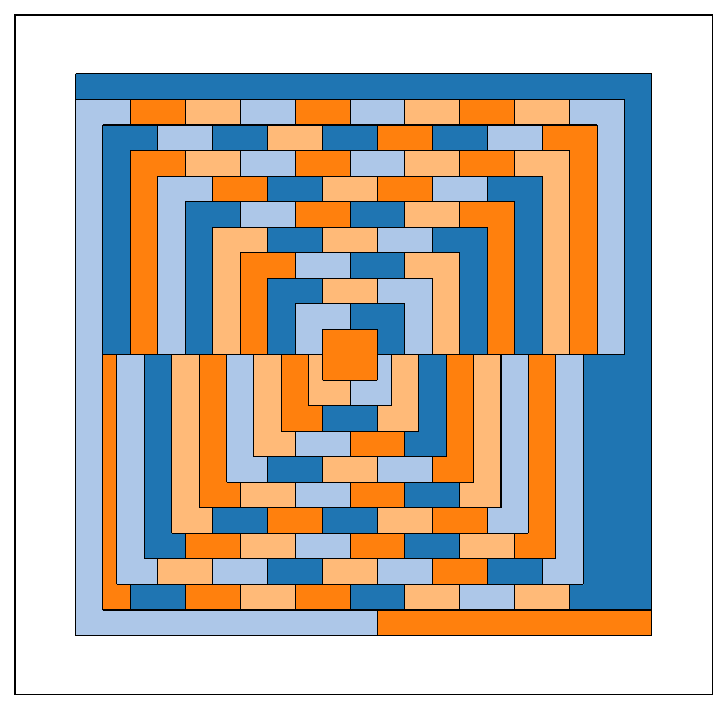}
			\includegraphics[trim=10 10 10 10, clip, scale=0.68]{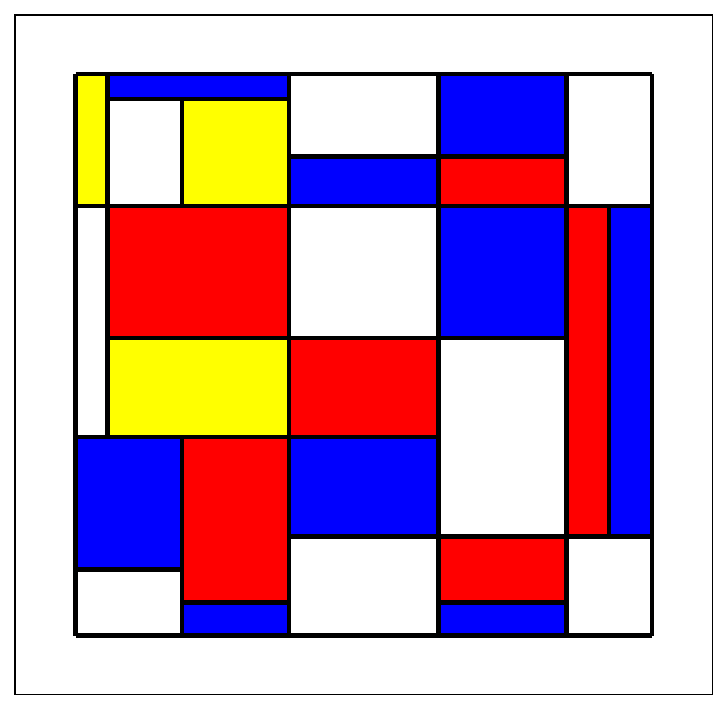}
		\end{center}
		\caption{The left image shows the planar embedding that was (jokingly) claimed to have a face chromatic number of five~\cite{Gardner1975}. A four colouring of this graph is shown. The right image shows a mock-up of Mondrian's \emph{Composition A} (1923).}
		\label{fig:joke}
	\end{figure}
	
	In 1975, Martin Gardner wrote a spoof article for \emph{Scientific American} to celebrate April Fool's Day~\cite{Gardner1975}. One of the false claims made in the article was that a graph had been discovered that featured a face chromatic number of five, therefore disproving the four colour theorem:
	
	\begin{quote}
		In November 1974 William McGregor, a graph theorist of Wappingers Falls, N.Y., constructed a map of 110 regions that cannot be colored with fewer than five colors. McGregor's technical report will appear in 1978 in the Journal of Combinatorial Theory, Series B.~\cite{Gardner1975}
	\end{quote}
		
	Of course, like all planar embeddings, a four colouring \emph{is} possible for this graph, as shown in Figure~\ref{fig:joke}.
	
	\section{Conclusions}
	
	The paper has reviewed and visualised several important results from the field of graph colouring. In general, optimally colouring a graph's nodes, edges, and faces is $\mathcal{NP}$-hard, though, as seen above, many special cases can be solved in polynomial time. For large, dense graphs, node and edge colourings can sometimes appear quite cluttered and difficult to interpret. However, careful node placement can improve this. On the other hand, face colourings are easier to visualise as they only apply to planar graphs.
	
	The four colour theorem, and the face colouring problem more generally, originated from the observation that, when colouring territories on a geographical map, no more than four colours are needed~\cite{Wilson2003}. Despite this, cartographers are not usually interested in limiting themselves to just four colours. Indeed, it is useful for maps to also satisfy other constraints, such as ensuring that all bodies of water (and no land areas) are coloured blue, and that disjoint areas of the same country (such as Alaska and the contiguous United States) receive the same colour. Such requirements can be modelled using the \emph{precolouring} and \emph{list colouring} problems, although they may well increase the required number of colours beyond four. Functionality for these problems has been included in the GCol library since version 2.2~\cite{Lewis2026GCol}.
	
	In 1890, Heawood~\cite{Heawood1890} provided a tight upper bound $H(h)$ on the number of colours needed to colour the faces of an $h$-holed torus for $h\geq 1$; specifically, $H(h) = \lfloor (7+\sqrt{1 + 48h}) / 2\rfloor$. It would be interesting to investigate the visualisation of such colourings. From an artistic perspective, it would also be interesting to explore the notion of optimal colourings in the cubist art of Picasso and Braque, or in the geometric art of Escher and Mondrian. An example of the latter is shown in Figure~\ref{fig:joke}.
	
	\bibliographystyle{plain}
	\bibliography{allrefs}

@article{Gardner1975, 
	year = {1975}, 
	volume = {232},
	number = {4},
	author = {Gardner, M.}, 
	title = {Mathematical Games: Six Sensational Discoveries that Somehow or Another have Escaped Public Attention}, 
	journal = {Scientific American}, 
	pages = {126-133}
}

@article{Dailey1980, 
	year = {1980}, 
	volume = {30}, 
	author = {Dailey, D.}, 
	title = {Uniqueness of Colorability and Colorability of Planar 4-Regular Graphs are {NP}-Complete}, 
	journal = {Combinatorial Theory}, 
	pages = {289-293}
}

@article{Smith2024, 
	doi = {10.5070/C64163843},
	year = {2024}, 
	volume = {4}, 
	number = {1}, 
	author = {Smith, D. and Myers, J. and Kaplan, C. and Goodman-Strauss, C.}, 
	title = {An Aperiodic Monotile}, 
	journal = {Combinatorial Theory} 
}

@article{Brinkmann2013,
  author    = {Brinkmann, B. and Coolsaet, K. and Goedgebeur, J. and M\'elot, H.},
  title     = {House of Graphs: a database of interesting graphs},
  journal   = {Discrete Applied Mathematics},
  volume    = {161},
  number    = {1-2},
  pages     = {311--314},
  year      = {2013},
  doi       = {10.1016/j.dam.2012.07.018},
  url       = {https://houseofgraphs.org/}
}

@MISC{Fischer2026,
  author = {Fischer, G.},
  title = {Image Triangulator, \url{https://snorpey.github.io/triangulation/}},
  year = {2026},
}

@book{diBattista1999,
  title={Graph Drawing: Algorithms for the Visualization of Graphs},
  author={di Battista, G. and Eades, P. and Tamassia, R. and Tollis, I.},
  year={1999},
  publisher={Prentice Hall}
}

@article{Fruchterman1991,
  title={Graph drawing by force-directed placement},
  author={Fruchterman, T. and Reingold, E.},
  journal={Software: Practice and Experience},
  volume={21},
  number={11},
  pages={1129--1164},
  year={1991}
}

@BOOK{Cranston2024,
  title = {Graph Coloring Methods},
  year = {2004},
  author = {Cranston, D.},
  address = {Richmond, Virginia},
  isbn = {979-218-46240-0},
  url = {\url{https://graphcoloringmethods.com/}}
}

@article{Lewis2025, 
	doi = {10.21105/joss.07871},
	url = {https://doi.org/10.21105/joss.07871}, 
	year = {2025}, 
	publisher = {The Open Journal}, 
	volume = {10}, 
	number = {108}, 
	pages = {7871}, 
	author = {Lewis, R. and Palmer, G.}, 
	title = {{GCol}: A High-Performance {P}ython Library for Graph Colouring}, 
	journal = {Journal of Open Source Software} 
}

@MISC{Lewis2026GCol,
  author = {Lewis, R.},
  title = {{GCol}: A Library for Graph Colouring v~2.2, \url{https://gcol.readthedocs.io/en/latest/}},
  year = {2026},
}

@INPROCEEDINGS{Johnson1974,
  author = {Johnson, D.},
  title = {Worst-case behavior of graph coloring algorithms},
  booktitle = {Proc. 5th Southeastern Conf. on Combinatorics, Graph Theory, and Computing, Utilitas Mathematicae},
  year = {1974},
  address = {Winnipeg, Canada},
  pages = {513-527}
}

@ARTICLE{Lewis2016,
  author = {Lewis, R. and Carroll, F.},
  title = {Creating Seating Plans: A Practical Application},
  journal = {Journal of the Operational Research Society},
  year = {2016},
  volume = {67},
  pages = {1353--1362},
  number = {11}
}

@ARTICLE{Aardel2002,
  author = {Aardel, K. and van~Hoesel, S. and Koster, A. and Mannino, C. and
	Sassano, A.},
  title = {Models and Solution Techniques for the Frequency Assignment Problems},
  journal = {4OR : Quarterly Journal of the Belgian, French and Italian Operations
	Research Societies},
  year = {2002},
  volume = {1},
  pages = {1--40},
  number = {4}
}

@ARTICLE{Appel1977,
  author = {Appel, K. and Haken, W.},
  title = {Solution of the Four Color Map Problem},
  journal = {Scientific American},
  year = {1977},
  volume = {4},
  pages = {108-121}
}

@ARTICLE{Brelaz1979,
  author = {Br\'{e}laz, D.},
  title = {New methods to color the vertices of a graph},
  journal = {Commun. ACM},
  year = {1979},
  volume = {22},
  pages = {251--256},
  number = {4}
}

@ARTICLE{Cambazard2012,
  author = {Cambazard, H. and Hebrard, E. and O'Sullivan, B. and Papadopoulos,
	A.},
  title = {Local search and constraint programming for the post enrolment-based
	timetabling problem},
  journal = {Annals of Operational Research},
  year = {2012},
  volume = {194},
  pages = {111-135}
}

@book{Lewis2021,
  author = {Lewis, R.},
  title = {A Guide to Graph Colouring: Algorithms and Applications},
  year = {2021},
  isbn = {978-3-030-81056-6},
  publisher = {Springer Cham},
  edition = {2nd}
}

@ARTICLE{Dupont2009,
  author = {Dupont, A. and Linhares, A. and Artigues, C. and Feillet, D. and
	Michelon, P. and Vasquez, M.},
  title = {The Dynamic Frequency Assignment Problem},
  journal = {European Journal of Operational Research},
  year = {2009},
  volume = {195},
  pages = {75-88}
}

@ARTICLE{Galinier1999,
  author = {Galinier, P. and Hao, J.},
  title = {Hybrid evolutionary algorithms for graph coloring},
  journal = {Journal of Combinatorial Optimization},
  year = {1999},
  volume = {3},
  pages = {379--397}
}

@book{Garey1979,
  author    = {Garey, M. and Johnson, D},
  title     = {Computers and Intractability: A Guide to the Theory of {NP}-Completeness},
  publisher = {W. H. Freeman and Company},
  address   = {San Francisco},
  year      = {1979},
  isbn      = {0-7167-1045-5}
}

@ARTICLE{Heawood1890,
  author = {Heawood, P.},
  title = {Map-Colour Theorems},
  journal = {Quarterly Journal of Mathematics},
  year = {1890},
  volume = {24},
  pages = {332-338}
}

@ARTICLE{Hoyler1981,
  author = {Holyer, I.},
  title = {The {NP}-completeness of edge-coloring},
  journal = {SIAM Journal on Computing},
  year = {1981},
  volume = {10},
  pages = {718-720}
}

@ARTICLE{Kirkman1847,
  author = {Kirkman, T.},
  title = {On a problem in combinations},
  journal = {Cambridge Dublin Math Journal},
  year = {1847},
  volume = {2},
  pages = {191-204}
}

@ARTICLE{Kubale1985,
  author = {Kubale, M. and Jackowski, B.},
  title = {A generalized implicit enumeration algorithm for graph coloring},
  journal = {Commun. ACM},
  year = {1985},
  volume = {28},
  pages = {412-418},
  number = {28}
}

@ARTICLE{Lewis2008,
  author = {Lewis, R.},
  title = {A Survey of Metaheuristic-based Techniques for University Timetabling
	Problems},
  journal = {OR Spectrum},
  year = {2008},
  volume = {30},
  pages = {167-190},
  number = {1}
}

@ARTICLE{Misra1992,
  author = {Misra, J. and Gries, D.},
  title = {A constructive proof of {V}izing's Theorem},
  journal = {Information Processing Letters},
  year = {1992},
  volume = {41},
  pages = {131-133}
}

@ARTICLE{Robertson1997,
  author = {Robertson, N. and Sanders, D. and Seymour, P. and Thomas, R.},
  title = {The Four Color Theorem},
  journal = {Journal of Combinatorial Theory, Series B},
  year = {1997},
  volume = {70},
  pages = {2-44}
}

@ARTICLE{deWerra1988,
  author = {de Werra, D.},
  title = {Some models of graphs for scheduling sports competitions},
  journal = {Discrete Applied Mathematics},
  year = {1988},
  volume = {21},
  pages = {47-65}
}

@BOOK{Wilson2003,
  title = {Four Colors Suffice: How the Map Problem Was Solved},
  publisher = {Penguin Books},
  year = {2003},
  author = {Wilson, R.}
}
	
\end{document}